\documentclass[11pt]{amsart}

\usepackage{amsfonts}
\usepackage[]{amsmath}
\usepackage[]{amssymb}
\usepackage[]{amsthm}
\usepackage[]{amsfonts}
\usepackage{ wasysym }
\usepackage{enumerate}

\newtheorem  {theorem}            {Theorem}[section]
\newtheorem  {lemma}[theorem]{Lemma}
\newtheorem  {definition}[theorem]{Definition}

\newtheorem {proposition}[theorem]{Proposition}

\newtheorem* {conjecture*}        {Conjecture}
\newtheorem* {theorem*}           {Theorem}

\theoremstyle{remark}

\newtheorem  {remark}[theorem]{Remark}
\newtheorem {remark*}[theorem]{Remark}
\newtheorem* {acknowledgements*}  {Acknowledgements}

\newcommand {\Ric}  {\operatorname{Rc}}
\newcommand {\R}    {\operatorname{\mathbb{R}}}

\newcommand {\mc}   {{H}}

\newcommand {\dist} {\operatorname{dist}}
\renewcommand {\div} {\operatorname {div}}

\newcommand {\tr} {\operatorname {tr}}

\renewcommand {\L} {{\mathcal{L}}}
\renewcommand {\H} {{\mathcal{H}}}
\renewcommand {\epsilon} {{\varepsilon}}

\newcommand {\Scal} {\operatorname {R}}
\newcommand {\C} {\operatorname {{C}}}

\newcommand {\graph} {\operatorname {graph}}

\newcommand{\la}{\langle}
\newcommand{\ra}{\rangle}

\title {Jenkins--Serrin--type results for the Jang equation}

\author{Michael Eichmair \& Jan Metzger}
\address{Michael Eichmair, Faculty of Mathematics, University of Vienna, Oskar-Morgenstern-Platz 1, 1090 Vienna, Austria }
\email{michael.eichmair@univie.ac.at}

\address{Jan Metzger, Institute of Mathematics, University of Potsdam,  Am Neuen Palais 10, 14469 Potsdam, Germany}
\email{jan.metzger@uni-potsdam.de}

\begin{document}

\maketitle

\begin{abstract} Let $(M, g, k)$ be an initial data set for the Einstein equations of general relativity. 

We prove that there exist solutions of the Plateau problem for marginally outer trapped surfaces (MOTSs) that are stable in the sense of MOTSs. This answers a question of G. Galloway and N. O'Murchadha raised in \cite{Galloway-OMurchadha:2008} and is an ingredient in the proof of the spacetime positive mass theorem \cite{spacetimePMT} given by L.-H. Huang, D. Lee, R. Schoen, and the first named author. 

We show that a \emph{canonical} solution of the Jang equation exists in the complement of the union of all weakly future outer trapped regions in the initial data set with respect to a given end, provided that this complement contains no weakly past outer trapped regions. The graph of this solution relates the area of the horizon to the global geometry of the initial data set in a non-trivial way.  

We prove the existence of a Scherk--type solution of the Jang equation outside the union of all weakly future or past outer trapped regions in the initial data set. This result is a natural exterior analogue for the Jang equation of the classical Jenkins--Serrin theory. 

We extend and complement existence theorems \cite{Jenkins-Serrin:1966, Jenkins-Serrin:1968, Spruck:1972, Nelli-Rosenberg:2002, Hauswirth-Rosenberg-Spruck:2009, Pinheiro:2009, Folha-Rosenberg:2012} for Scherk--type  constant mean curvature graphs over polygonal domains in $(M, g)$, where $(M, g)$ is a complete Riemannian surface. We can dispense with the a priori  assumptions that a sub solution exists and that $(M, g)$ has particular symmetries. Also, our method generalizes to higher dimensions.  
\end{abstract}


\section{Introduction and Notation} \label{sec:intro}

In this paper we prove a series of results concerning the geometric theory of the Jang equation. We apply our insights to obtain an optimal extension to general Riemannian surfaces of the classical Jenkins-Serrin-Spruck theory \cite{Jenkins-Serrin:1968, Spruck:1972} on the existence of Scherk-type minimal and constant mean curvature graphs. In fact, our methods allow us to remove two conditions from J. Spruck's pioneering result \cite{Spruck:1972} in $\R^2$.  \\

In Section \ref{sec:Plateau} we show that there exist stable solutions of the Plateau problem for marginally outer trapped surfaces. This answers a question of G. Galloway and N. O'Murchadha raised in \cite{Galloway-OMurchadha:2008}. This result is subtle in view of the non-variational nature of these surfaces. The proof is based on the fact that using the existence theory in \cite{Eichmair:2009-Plateau}, we can construct \emph{ordered} families of solutions of the Plateau problem. This result is applied in the recent proof of the spacetime positive mass theorem by L.-H. Huang, D. Lee, R. Schoen, and the first named author in \cite{spacetimePMT}. Stability in the existence theory for closed MOTSs was concluded in \cite[Section 4]{Andersson-Metzger:2009} (see also \cite[Section 3.6]{Andersson-Eichmair-Metzger:2010} for a simplification) by a different argument that does not extend to the present situation.  \\

In Sections \ref{sec:scherk}, \ref{sec:canonical}, and \ref{sec:uniqueness} we develop the geometric theory of the Jang equation pioneered by R. Schoen and S.-T. Yau in \cite{Schoen-Yau:1981-pmt2} to prove the existence of non-trivial and, in some cases, canonical Scherk-type solutions of the Jang equation in the complement of the total weakly future outer trapped region and the total weakly past outer trapped region. \\

In Section \ref{sec:jss} we employ techniques from the geometric theory of the Jang equation, in particular the capillarity regularization and  the geometric blow up analysis from \cite{Schoen-Yau:1981-pmt2} and ideas from the solution of the non-variational Plateau problem for marginally outer trapped surfaces  in \cite{Eichmair:2009-Plateau}, to the classical Jenkins--Serrin--Spruck problem \cite{Jenkins-Serrin:1966, Jenkins-Serrin:1968, Spruck:1972} of finding necessary and sufficient conditions for a domain in a Riemannian surface to support a Scherk-type constant mean curvature graph. In the case of positive mean curvature, we are able to dispense with the a priori assumption that the domain admit a sub solution which is required in the foundational paper by J. Spruck \cite{Spruck:1972} (in $\R^2$) and its recent extension to domains in $\mathbb{S}^2$ and $\mathbb{H}^2$ by L. Hauswirth, H. Rosenberg, and J. Spruck \cite{Hauswirth-Rosenberg-Spruck:2009}, and to domains in Hadamard manifolds \cite{Folha-Rosenberg:2012} by A. Folha and H. Rosenberg. Moreover, our results are valid in arbitrary complete Riemannian surfaces. Thus the existence of a Scherk-type constant mean curvature graph above a Riemannian surface is fully reduced to a (generically) finite set of inequalities relating area and circumference of certain polygons that can be inscribed into the domain: the Jenkins-Serrin-Spruck flux conditions. Our approach here does not distinguish between minimal and (positive) constant mean curvature graphs. In the case of minimal graphs, we recover the results by B. Nelli and H. Rosenberg \cite{Nelli-Rosenberg:2002} (in $\mathbb{H}^2$) and by A. Pinheiro \cite{Pinheiro:2009} (in general Riemannian surfaces). Our methods carry over to higher dimensions. A more detailed overview of the literature and a precise statement of our result are given in Section \ref{sec:jss}. \\

In the appendices we collect several results that are used in this paper. In Appendix \ref{sec:ige} we deduce the interior gradient estimate for solutions of the prescribed mean curvature equation in low dimensions from the regularity theory for almost minimal boundaries. In Appendix \ref{sec:boundary} we characterize those domains that support infinite boundary value solutions of the prescribed mean curvature equation \emph{without} making an a priori assumption on the regularity of the boundary. In Appendix \ref{sec:equicontinuous} we observe a simple and useful consequence of the classical compactness and regularity theory for almost minimal boundaries: the horizontal parts of the unit normal vector fields of solutions to the prescribed mean curvature equation are equicontinuous in low dimensions. \\ 

We proceed by introducing some notation and conventions.   \\

Let $(M, g)$ be a connected Riemannian manifold of dimension $n$, with $3 \leq n \leq 7$, and let $k$ be a symmetric $(0, 2)$-tensor on $M$. In the context of the Cauchy problem for the Einstein equations in general relativity, $k$ is referred to as the (spacetime) second fundamental form tensor and the triple $(M, g, k)$ is called an initial data set. We often require that $(M, g, k)$ is asymptotically flat, i.e. that the complement of some compact subset of $M$ consists of finitely many connected components $N_1, \ldots, N_m$, called the ends, each one diffeomorphic to $\R^n \setminus \overline{ B_1(0) }$ and such that in the corresponding coordinate systems the metric tensor $g_{ij}$ converges to the Euclidean metric $\delta_{ij}$ and the second fundamental form tensor $k_{ij}$ to zero. More precisely, we require that
\begin{align*} \label{eqn:AF}
|g_{ij} - \delta_{ij}| + |x||\partial_k g_{ij}| = O(|x|^{-q})
\end{align*}
and 
\begin{align*}
k_{ij} = O(|x|^{-q-1})
\end{align*}
as $|x| \to \infty$ for some 
\[
q > \frac{n-2}{2}.
\]
When $n=3$, we ask in addition that for some $\beta >2$, 
\begin{align*} \label{eqn:AF3}
\tr_g(k) = g^{ij} k_{ij} = O(|x|^{-\beta})
\end{align*}
as $x \to \infty$. This last condition is imposed so that certain barriers for the Jang equation can be constructed far out in the asymptotically flat ends, cf. Appendix \ref{sec:barriers}. \\

Let $(M, g, k)$ be an initial data set, and let $\Sigma \subset M$ be a two-sided hypersurface with unit normal vector field $\nu$. The future outer and past outer expansion scalars of $\Sigma$ are defined, respectively, as $\theta^+_\Sigma = H_\Sigma + \tr_\Sigma (k)$ and $\theta^-_\Sigma = H_\Sigma - \tr_\Sigma (k)$. Here, $H_\Sigma$ is the mean curvature scalar of $\Sigma$ computed as the tangential divergence of $\nu$ and $\tr_\Sigma (k)$ is the trace of $k$ restricted to the tangent space of $\Sigma$. The hypersurface $\Sigma$ is called future outer trapped (respectively past outer trapped) if $\theta^+_\Sigma <0$ ($\theta^-_\Sigma <0$) everywhere on $\Sigma$, weakly future outer trapped (respectively weakly past outer trapped) if $\theta^+_\Sigma \leq 0$ ($\theta^-_\Sigma \leq 0$) everywhere on $\Sigma$, and is called a marginally future outer trapped surface or MOTS for short (respectively marginally past outer trapped surface or MITS for short) if $\theta^+_\Sigma = 0$ ($\theta^-_\Sigma =0$) everywhere on $\Sigma$. Except for Section \ref{sec:Plateau}, the MOTSs and MITSs appearing in this paper will be closed and in fact boundaries of sets $\Omega$ that contain part of an end $\{x \in N_i : |x| \geq r_0\}$ for some $r_0 \geq 1$ and $i \in \{1, \ldots, m\}$. If we write a hypersurface $\Sigma$ as the (relative) boundary of a set $\Omega$, say $\Sigma = \partial \Omega$ in $U$ where $U$ is an open subset of $M$, we will use the unit normal field $\nu$ of $\Sigma$ pointing into $\Omega$ to compute the scalar mean curvature $H_\Sigma$ and the expansion scalars $\theta^\pm_\Sigma$ unless otherwise noted. \\
 
Let $(M, g, k)$ be a complete asymptotically flat manifold of dimension $3 \leq n \leq 7$. Fix one of the ends, say $N_1$. It is easy to see that $H_{S_r} > |\tr_{S_r} k|$ for all $r \geq r_0$, provided that $r_0 \geq 1$ is sufficiently large. Here, $S_r = \{x \in N_1 : |x| = r\}$ is the coordinate sphere of radius $r$ in $N_1$. The mean curvature scalar $H_{S_r}$ is computed as the tangential divergence of the unit normal pointing into the end. Let $M_- \subset M$ (respectively $M_+ \subset M$) be the interior of the intersection of all open subsets $\Omega \subset M$ that contain $\{x \in N_1 : |x| \geq r_0\}$ and which have smooth compact embedded boundary satisfying $H_{\partial \Omega} + \tr_{\partial \Omega} (k) \leq 0$ ($H_{\partial \Omega} - \tr_{\partial \Omega} (k) \leq 0$). It follows from \cite{Andersson-Metzger:2009} in dimension $n=3$ and from \cite{Eichmair:2010} in dimensions $3 \leq n \leq 7$ that $M_-$ (respectively $M_+$) has smooth compact embedded boundary such that $H_{\partial M_-} + \tr_{\partial M_-} (k) = 0$ ($H_{\partial M_+} -\tr_{\partial M_+} (k) = 0$). Thus $\partial M_-$ is a MOTS and $\partial M_+$ is a MITS. The complements of the regions $M_-$ and $M_+$ are the total weakly future outer trapped region and the total weakly past outer trapped region of $(M, g, k)$ with respect to the chosen end, respectively. If $(M, g, k)$ has more than one end, then both $\partial M_-$ and $\partial M_+$ are non-empty and each of them separates the portion $\{x \in N_1 : |x| \geq r_0\}$ of $N_1$ from $\bigcup_{i=2}^m \{x \in N_i : |x| \geq r_0\}$ provided $r_0$ is sufficiently large. \\

Given an open subset $\Omega \subset M$ and a function $u \in \C^2(\Omega)$ we define, in local coordinates near a point $x \in \Omega$, 
\begin{align*} 
H(u) = \left( g^{ij} - \frac{u^i u^j}{1 + |D u|^2}\right) \frac{D^2_{ij} u }{\sqrt{1 + |D u|^2}}
\end{align*}
and 
\begin{align*} 
\tr (k) (u)=  \left( g^{ij} - \frac{u^i u^j}{1 + |D u|^2}\right) k_{ij}. 
\end{align*}  
Here, all geometric operations (raising indices, gradient, length of gradient, Hessian) are with respect to $g$. These definitions are independent of the particular coordinate system used. The function $H(u)$ at $x \in \Omega$ is the scalar mean curvature of the graph $G = \{(x, u(x)) : x \in \Omega\}$ of $u$ in the Riemannian product $(\Omega \times \R, g + dx^{n+1} \otimes dx^{n+1})$ at the point $(x, u(x))$ with respect to the downward pointing unit normal, and $\tr(k)(u)$ evaluated at $x$ is the trace of $k$ (extended to $\Omega \times \R$ by zero in the vertical direction) over the tangent space of this graph at $(x, u(x))$. If
\begin{eqnarray} \label{Jang} 
H(u) + \tr(k)(u) = 0,
\end{eqnarray}
then $G$, with its downward orientation, is a MOTS in the new initial data set $(M \times \R, g + d x^{n+1} \otimes d x^{n+1}, k)$. Equation (\ref{Jang}) is known as the Jang equation. \\

For background material on MOTSs, MITSs, and the Jang equation we refer the reader to the survey article \cite{Andersson-Eichmair-Metzger:2010}.\\ 

\section*{Acknowledgments} We would like to thank L. Andersson, G. Galloway, G. Huisken, and M. Mars for their interest in this work and for valuable discussions. M. Eichmair would like to thank R. Beig, P. Chru\'sciel, J. Grant, M. Heinzle, and W. Simon of the group for Gravitational Physics at the University of Vienna, as well as the wonderful Erwin 
Schr\"odinger Institute, for kindly providing pleasant and stimulating working conditions for him in the summer of 2011. He also gratefully acknowledges the support of a Clay Liftoff Fellowship in the summer of 2008, when some of the results in this paper were first conceived, and the support of NSF grant DMS-0906038 and of SNF grant 200021-140467. The second named author gratefully acknowledges support by the DFG grant ME3816/1-1.


\section{Stability of solutions of the Plateau problem} \label{sec:Plateau}

The notion of stability for MOTSs with boundary in the definition below is natural in view of the notion of stability for closed MOTSs that has been
introduced and studied systematically in \cite{Andersson-Mars-Simon:2008}.

\begin{definition} [Cf. \protect{\cite[Section 2]{Galloway-OMurchadha:2008}}]
  \label{def:stability_of_MOTS}
Let $(M, g, k)$ be an initial data set and consider a two-sided hypersurface
$\Sigma \subset M$ with boundary $\partial \Sigma = \Gamma$ and with a designated
``outward" unit normal vector field $\nu$. Assume that $\Sigma$ is a MOTS.  Then $\Sigma$ is said
to be {\it stable in the sense of MOTSs} if there exists a smooth function $f$
that is positive in the interior of $\Sigma$, vanishes on the boundary
$\Gamma$, and is such that $L_\Sigma f \geq 0$. Here,
  \begin{equation*}
    L_\Sigma \phi = - \Delta_\Sigma \phi + 2 \la X,
    D_\Sigma \phi \ra + \left( \tfrac{1}{2} \Scal_\Sigma - \tfrac{1}{2} |h +
      k|_\Sigma^2 - J (\nu) - \mu + \div_\Sigma X - |X|^2 \right) \phi,
  \end{equation*}
  where $X$ is the tangential part of the vector field dual to $k(\nu, \cdot)$ on
$\Sigma$, $h$ is the second fundamental form of $\Sigma$ (with trace
$\mc_\Sigma$), $D_\Sigma \phi$ is the tangential gradient of $\phi$ along
$\Sigma$, $\Scal_\Sigma$ is  the scalar curvature of $\Sigma$, $\mu = \frac{1}{2}
\left (\Scal_M + (\tr_M (k))^2 - |k|_M^2\right)$ is the local mass density, and
where $J = \div(k - \tr_M(k)g)$ is the local current density. Equivalently, $\Sigma$ is stable in the sense of MOTSs if, and only if, the principal eigenvalue of $L_\Sigma$ is non-negative. 
\end{definition}

For a careful discussion of principal eigenvalues of (not necessarily self-adjoint) elliptic operators, we refer the reader to \cite[Sections 3.6 and 3.7]{Pinsky:1995}. \\

We recall the following existence theorem for MOTSs spanning a given boundary:

\begin{theorem} [\cite{Eichmair:2009-Plateau}] \label{thm:existencePlateau}
Let $(M, g, k)$ be a complete initial data set of dimension $n$ with $2 \leq n \leq 7$. Let $\Omega \subsetneq M$ be a bounded open set with smooth boundary $ \partial \Omega$. Let $\Gamma \subset \partial \Omega$ be a non-empty smooth closed embedded submanifold of $\partial \Omega$ such that $\partial \Omega \setminus \Gamma = \partial_- \Omega \dot \cup \partial_+ \Omega$ for disjoint non-empty relatively open subsets $\partial_- \Omega, \partial_+ \Omega$ of $\partial \Omega$. Assume that $\mc_{\partial \Omega} + \tr_{\partial \Omega} k < 0$ near $\partial_- \Omega$ with the mean curvature computed as the tangential divergence pointing into $\Omega$ and that $\mc_{\partial \Omega} + \tr_{\partial \Omega} k > 0$ near $\partial_+ \Omega$ with the mean curvature scalar computed as the tangential divergence of the unit normal pointing out of $\Omega$. Then there exists a smooth hypersurface $\Sigma \subset \Omega$ with boundary $\Gamma$ that is an almost minimizing relative boundary in $\Omega$ and such that $\Sigma$ is a MOTS with respect to the unit normal pointing towards $\partial_+ \Omega$.
\end{theorem}

The following theorem answers a question posed in \cite[Section 3]{Galloway-OMurchadha:2008}. It is an ingredient in the proof of the spacetime positive mass theorem given in \cite{spacetimePMT}. 

\begin{theorem}
Assumptions as in Theorem \ref{thm:existencePlateau}. Then there exists a solution $\Sigma$ of the Plateau problem for MOTSs in $\Omega$ with boundary $\Gamma$ that is stable in the sense of MOTSs.
\begin{proof} Given $\epsilon >0$ small, let $\Gamma^\epsilon = \{ \theta \in \partial_+ \Omega: \dist_{\partial \Omega} (\theta, \Gamma) = \epsilon\}$. It follows from the construction in the proof of Theorem \ref{thm:existencePlateau} in \cite[Chapter 4]{Eichmair:2009-Plateau} that the MOTSs $\Sigma^\epsilon \subset \Omega$ spanning $\Gamma^\epsilon$ are (strictly) ordered. To see this, recall that the open subset of $\Omega$ whose relative boundary is $\Sigma^\epsilon$ is the geometric limit as $t \searrow 0$ of downward translations of the regions lying above the graphs $u_{t}^\epsilon \in \C^\infty_{loc} (\Omega)$ constructed in \cite[Lemma 4.2]{Eichmair:2009-Plateau}. Given $t >0$ and $0 < \epsilon < \epsilon'$ small, we have that $\mathcal{S}_{\overline u^{\epsilon'}_{t}} \subset \mathcal{S}_{\overline u^{\epsilon}_{t}}$ (in the notation of \cite{Eichmair:2009-Plateau}) and hence $u^{\epsilon'}_{t} \leq  u^\epsilon_{t}$. It follows that the regions above the graphs are ordered so that $\Sigma^{\epsilon'}$ lies to one side (towards $\partial_+ \Omega$) of $\Sigma^{\epsilon}$. The geometric maximum principle shows that components of $\Sigma^{\epsilon'}$ and $\Sigma^{\epsilon}$ that span components of $\Gamma^\epsilon$ and $\Gamma^{\epsilon'}$ cannot touch unless they coincide. We will discard all extraneous closed components of $\Sigma^\epsilon$. This does not change that each $\Sigma^\epsilon$ is a relative boundary, nor that the $\Sigma^\epsilon$'s are ordered.  

The geometric compactness properties of the almost minimizing relative boundaries $\Sigma^\epsilon$ show that as $\epsilon \searrow 0$, the $\Sigma^\epsilon$ converge smoothly and with multiplicity one to an embedded MOTS $\Sigma$ that spans $\Gamma$. We claim that this MOTS $\Sigma$ is stable in the sense of MOTSs. To see this, let $U, V, W \subset \Sigma$ be non-empty open subsets with smooth boundaries such that $U \Subset V \Subset W \Subset \text{int } \Sigma$. Let $\nu$ be the unit normal vector field of $\Sigma$ that points towards $\partial_+ \Omega$. (This makes sense because $\Sigma$ is a relative boundary in $\Omega$ spanning $\Gamma$.) By assumption, there exist positive functions $f^\epsilon \in \C^\infty (W)$ for $\epsilon >0$ small with 
$\{\exp_\theta \left( f^\epsilon (\theta) \nu (\theta) \right) : \theta \in V \} \subset \Sigma^\epsilon$
and such that $f^\epsilon \to 0$ with all derivatives on compact subsets of $W$. Because $\Sigma$ and $\Sigma^\epsilon$ both satisfy the MOTS equation, $f^\epsilon$ is solution of a homogeneous linear elliptic equation in divergence form on $V$. The operator describing the linearization of the equation at the function that vanishes identically is $L_\Sigma$. Arguing exactly as in \cite[p. 333]{Simon:1987}, using Harnack theory, it follows that the functions $f^\epsilon$ can be rescaled (so their infimum is one on $V$, say) so as to converge smoothly to a positive function $f_V \in C^\infty(V)$ with $L_\Sigma (f_V) = 0$. This implies that $U$ is stable in the sense of MOTSs. To see this, let $\lambda$ be the first principal eigenvalue of $L_\Sigma|_U$ and let $h \in \C^\infty (\bar U)$ be the corresponding first (Dirichlet) eigenfunction so that $L_\Sigma h = \lambda h$. We recall that the first principal eigenvalue is simple and that the corresponding eigenfunctions do not change signs. By scaling, using that $g$ vanishes on the boundary of $U$ and that $f_V$ is positive on $V$ and that $U \Subset V$, we may assume that $0 < h \leq f_V$ on $U$ with equality at some point. The maximum principle then implies that $\lambda \geq 0$. We conclude that every open subset $U \Subset \text{int } \Sigma$ is stable in the sense of MOTSs. Using that the principal eigenvalue of an elliptic operator depends continuously on the operator and the domain, it follows that $\Sigma$ is stable in the sense of MOTSs. 
\end{proof}
\end{theorem}


\section{Scherk-type solutions of the Jang equation} \label{sec:scherk}

The content of the following proposition is similar to that of Theorem 3.1 in \cite{Metzger:2010-blowup}. We include an alternative proof here as preparation for the more general and difficult Theorem \ref{thm:JangScherk} below. The modification of the data near the boundary in our proof is much less delicate than that in \cite{Metzger:2010-blowup}. The regions $M_-$ and $M_+$ in the statement of Proposition \ref{prop:existenceblowup} and Theorem \ref{thm:JangScherk} below are defined in the introduction.  

\begin{proposition} [Cf. Theorem 3.1 in \cite{Metzger:2010-blowup} when $n=3$] \label{prop:existenceblowup} Let $(M, g, k)$ be a complete asymptotically flat initial data set of dimension $n$,  $3 \leq n \leq 7$, and assume that $\partial M_-$ and $\partial M_+$ are disjoint. There exists a smooth solution $u : M_- \cap M_+ \to \mathbb {R}$ of the Jang equation $\mc(u) + \tr(k)(u) = 0$ such that $u(x) \to 0$ as $x \to \infty$ in the  end of $(M, g, k)$ contained in $M_- \cap M_+$ and such that $u (x) \to \pm \infty$ as $\dist(x, \partial M_\pm) \to 0$. 

\begin{proof} We abbreviate $\Omega = M_- \cap M_+$. Let $\chi \in \C^\infty_c (M)$ be such that $\chi \equiv \pm 1$ near $\partial M_\pm$. Given $\epsilon \in (0, 1)$ we define $k^\epsilon = k + \epsilon \chi$. Note that $H_{\partial M_+} + \tr_{\partial M_+}(k^\epsilon) > 0$ and $H_{\partial M_-} + \tr_{\partial M_-} (k^\epsilon) < 0$. 

Arguing exactly as in \cite[Chapter 3]{Eichmair:2009-Plateau} and \cite[Chapters 3 and 4]{Eichmair:2010}, using also the argument in Appendix \ref{sec:barriers} and Footnote \ref{footnote:noncompactdomains}, we see that for every $\epsilon \in (0, 1)$ there exists a connected open subset $\Omega_0^\epsilon \subset \Omega$ containing the chosen end of $(M, g, k)$ and a solution $u^\epsilon \in \C^\infty_{loc} (\Omega_0^\epsilon)$ of the Jang equation $H(u^\epsilon) + \tr(k^\epsilon) (u^\epsilon) = 0$ on $\Omega_0^\epsilon$ with the following properties: 
\begin{enumerate} [(i)]
\item \label{prop:topstab} We have that $\{x \in \Omega : |x| > \Lambda\} \subset \Omega_0^\epsilon$ for some $\Lambda \geq 1$ that is independent of $\epsilon \in (0, 1)$. The topological boundary $\partial \Omega_0^\epsilon$ of $\Omega_0^\epsilon$ is a smooth properly embedded hypersurface in $\Omega$ whose components are either marginally inner trapped or marginally outer trapped with respect to the unit normal pointing into $\Omega_0^\epsilon$. The components are $\lambda$--minimizing with $\lambda = 1 + 2 n  \sup_{x \in \Omega, \epsilon \in (0, 1) } |k^\epsilon (x)|$ in $\Omega$ (in the language of \cite{Duzaar-Steffen:1993a}) and stable (in the sense of (\ref{eqn:almoststability})).
\item We have that $u(x) \to 0$ as $x \to \infty$ in $\Omega$. We have that $u^\epsilon(x)$ diverges to plus infinity if $x \in \Omega_0^\epsilon$ approaches a marginally inner trapped component of the boundary of $\Omega_0^\epsilon$, and to minus infinity if $x \in \Omega$ converges to a marginally outer trapped boundary component.
\item \label{prop:stablegraph} The graphs $\{(x, u^\epsilon (x)) : x \in \Omega_0^\epsilon \}$ of $u^\epsilon$ are complete hypersurfaces of $M \times \R$ that are stable and  $\lambda$--minimizing in $\Omega \times \R$. 
\end{enumerate}
The $\lambda$--minimizing property and stability of the graphs asserted in (\ref{prop:stablegraph}) implies that of the components of $\partial \Omega_0^\epsilon$ in (\ref{prop:topstab}), cf. \cite[p. 254]{Schoen-Yau:1981-pmt2}, \cite[Lemma A.1] {Eichmair:2009-Plateau}, and the discussion in Appendix \ref{sec:boundary}. 

As in \cite{Eichmair:2010} or the discussion in Appendix \ref{sec:boundary}, we see that the stability and the almost minimizing property (via uniform local mass bounds) lead to curvature estimates for these graphs that are independent of $\epsilon \in (0, 1)$. We now let $\epsilon \searrow 0$ and pass the graphs of $u^\epsilon$ to a smooth, properly embedded\footnote{The almost minimizing property by itself does not lead to uniform curvature estimates near $\partial \Omega$. As in Appendix \ref{sec:boundary}, we use the stability and the completeness for that. Once we know that a smooth limit exists, we can use the almost minimizing property to rule out sheeting.} subsequential limit that contains a connected complete graphical component whose domain contains the asymptotically flat end and which satisfies all the above properties with $\epsilon = 0$. Using the mean value theorem and that there are no MOTSs or MITSs (with respect to $k$) in $\bar \Omega$ besides $\partial M_-$ and $\partial M_+$ we see that this graphical component has all the properties asserted in the conclusion of the theorem. 
\end{proof}
\end{proposition} 

\begin{remark} Both in Proposition \ref{prop:existenceblowup} above and in Theorem \ref{thm:JangScherk} below, the outermost property of $M_-$ and $M_+$ prevents the domain of the sought-after graphical solution of the Jang equation from ``popping outward". If we think of the results in this section in analogy with the classical Jenkins--Serrin theory \cite{Jenkins-Serrin:1966, Jenkins-Serrin:1968}, then this outermost property takes the place of the Jenkins--Serrin flux conditions (\ref{eqn:totalflux}), (\ref{eqn:partialfluxA}), and (\ref{eqn:partialfluxB}) in the statement of Theorem \ref{thm:mainjss} below (with $H_0=0$). 
\end{remark}

\begin{theorem} \label{thm:JangScherk} Let $(M, g, k)$ be a complete asymptotically flat initial data set of dimension $n$ where $3 \leq n \leq 7$. Assume that $\partial M_+$ and $\partial M_-$ are both non-empty and that they intersect transversely. There exists a smooth solution $u : M_- \cap M_+ \to \R$ of the Jang equation 
$H(u) + \tr(k) (u) = 0$ such that $u(x) \to 0$ as $x \to \infty$ in the  end of $(M, g, k)$ contained in $M_- \cap M_+$ and such that $u(x) \to \pm \infty$ as $x \in M_- \cap M_+$ approaches an interior point $y$ of a component of  $M_{\mp} \cap \partial M_{\pm} \subset \partial (M_- \cap M_+)$. The topological closure of the graph $\{ (x, u(x)) : x \in M_- \cap M_+ \}$ of $u$ in $M \times \R$ is a smooth properly embedded hypersurface with manifold boundary $(\partial M_+ \cap \partial M_-) \times \R$.    
\begin{proof} Let $\Omega = M_- \cap M_+$.
We denote by $\nu$ the unit normal vector field of the hypersurface $M_{\mp} \cap \partial M_\pm$ pointing into $\Omega$. Let $ \partial_- \Omega = M_+ \cap \partial M_-$ and $\partial_+ \Omega = M_- \cap \partial M_+$. Note that $\partial_- \Omega$ and $\partial_+\Omega$ are manifolds with boundary and that their manifold boundaries coincide. The topological boundary of $\Omega$ is the disjoint union of $\partial_- \Omega$, $\partial_+ \Omega$, and $\partial M_- \cap \partial M_+$. Given a component $\gamma$ of  $\partial_\pm \Omega$, let $\Theta_\gamma \in \C^\infty(\bar \gamma)$ be positive in the interior of $\gamma$ and zero on its manifold boundary. 

Let $\epsilon \in (0, 1)$ be so small that the sets $\{ \exp_{y}  t \Theta_\gamma (y) \nu(y) : y \in \text{int } \gamma \text { and } t \in (0, 2\epsilon) \} \subset  \Omega$ are disjoint as $\gamma$ ranges over the components of $\partial_\pm \Omega$. Let $\chi^\epsilon \in \C^\infty_{loc} (\Omega)$ with values in $[-1, 1]$ be supported in the union of all these sets and such that 
$
\chi^\epsilon (y) \equiv \pm 1$ on 
\[
\text{Cr}_{\gamma}^\epsilon = \{ \exp_{y}  t \Theta_\gamma (y) \nu(y) : y \in \text{int } \gamma \text { and } t \in (0, \epsilon)\}
\] when $\gamma$ is a component of $\partial_\pm \Omega$. 

Let $\gamma$ be a component of $\partial_- \Omega$. Consider the hypersurface $$\{ (\exp_y \epsilon ( 1 - e^{-h})  \Theta(y)  \nu (y), h) : y \in \text{int } \gamma \text{ and } h \in (0, \infty)\} \subset \Omega \times \R.$$ It has piecewise smooth manifold boundary consisting of $\gamma \times \{0\}$ and $(\partial \gamma) \times [0, \infty)$. It is the graph of a positive function $\overline u_{\gamma}^\epsilon \in \C^\infty_{loc} (\text{Cr}_{\gamma}^\epsilon)$. We have that $H(\overline u_{\gamma}^\epsilon) + \tr(k) (\overline u_{\gamma}^\epsilon) = O(\epsilon)$.    

Similarly, let $\gamma$ be a component of $\partial_+ \Omega$. Consider the hypersurface $$\{ (\exp_y \epsilon ( 1 - e^{h})  \Theta(y)  \nu (y), h) : y \in \text{int } \gamma \text{ and } h \in (- \infty, 0)\} \subset \Omega \times \R.$$ Its topological boundary is the union of $\gamma \times \{0\}$ and $(\partial \gamma) \times (- \infty, 0]$. This hypersurface is the graph of a negative function $\underline u_{\gamma}^\epsilon \in \C^\infty_{loc} (\text{Cr}_{\gamma}^\epsilon)$. We have that $H(\underline u_{\gamma}^\epsilon) + \tr(k) (\underline u_{\gamma}^\epsilon) = O(\epsilon)$.

Choose $c >0$ so that for all $\epsilon >0$ sufficiently small $H(\overline  u_{\gamma}^\epsilon) + \tr(k^\epsilon) (\overline u_{\gamma}^\epsilon) < - 2\epsilon$ on $\text{Cr}_{\gamma}^\epsilon$ when $\gamma$ is a component of $\partial_-\Omega$ and such that $H(\underline  u_{\gamma}^\epsilon) + \tr(k^\epsilon) (\underline u_{\gamma}^\epsilon) >  2\epsilon$ on $\text{Cr}_{\gamma}^\epsilon$ when $\gamma$ is a component of $\partial_+ \Omega$. Here, $k^\epsilon = k + c \epsilon \chi^\epsilon$.

Given $t>0$, the functions $- \frac{\epsilon}{t} + \overline u_{\gamma}^\epsilon$ and $\frac{\epsilon}{t} + \underline u_{\gamma}^\epsilon$ are, respectively, super and sub solutions of the regularized Jang equation $H(u) + \tr(k^\epsilon)(u) = t \,  u$ on $\text{Cr}_{\gamma}^\epsilon$. Let $C >0$ be a constant greater than $n \sup_{\epsilon \in (0, 1), x \in \Omega} |k^\epsilon(x)|$. Then $\frac{C}{t}$ and $- \frac{C}{t}$ are constant super and sub solutions of this equation. 

As in \cite[Chapters 3 or 4]{Eichmair:2009-Plateau} one sees that for every $t>0$ there exists a (Perron) solution $u_{t}^\epsilon \in \C^\infty _{loc} (\Omega)$ of $H(u_{t}^\epsilon) + \tr(k^{\epsilon}) (u_{t}^\epsilon) = t \, u_{t}^\epsilon$ such that $- \frac{C}{t} \leq u_{t}^\epsilon \leq \frac{C}{t}$ on $\Omega$, such that $u_{t}^\epsilon \leq - \frac{\epsilon}{t} + \overline u_{\gamma}^\epsilon $ on $\text{Cr}_{\gamma}^\epsilon$ for components $\gamma$ of $\partial_-\Omega$, such that $\frac{\epsilon}{t} + \underline u_{\gamma}^\epsilon \leq u_{t}^\epsilon$ on $\text{Cr}_{\gamma}^\epsilon$ for components $\gamma$ of $\partial_+ \Omega$, and such that $|u_t^\epsilon (x)| \leq b_\Lambda(|x|)$ on $\{x \in \Omega : |x| > \Lambda\}$ for some $\Lambda \geq 1$ sufficiently large. Here, $b_\Lambda$ is as in Appendix \ref{sec:barriers}.\footnote{\label{footnote:noncompactdomains} The results in \cite{Eichmair:2009-Plateau} are for compact domains. One way to obtain $u_t^\epsilon$ that decays to zero in the asymptotically flat end is as a limit of solutions $u_t^{\epsilon, r}$ on $\{x \in \Omega : |x| \leq r\}$ with zero boundary values on $\{x \in \Omega : |x| = r\}$ as $1 \leq r \to \infty$. Note that $|u_{t}^{\epsilon, r}(x)| \leq b_\Lambda (|x|)$ on $\{x \in \Omega : \Lambda < |x| \leq r\}$ for a fixed $\Lambda \geq 1$ sufficiently large.} 

It follows that if $y$ is an interior point of $\partial_-\Omega$, then $\limsup_{x \to y, x \in \Omega} u_{t}^\epsilon (x) \leq - \frac{\epsilon}{t}$ and that if $y$ is an interior point of $\partial_+ \Omega$, then $\liminf_{x \to y, x \in \Omega} u_{t}^\epsilon (x) \geq \frac{\epsilon}{t}$.

The mean curvature of the graphs of $u_{t}^\epsilon$ is bounded by $2C$ so that they are $2C$--minimizing (in the language of \cite{Duzaar-Steffen:1993a}) in $\Omega \times \R$, cf. \cite[Appendix A]{Eichmair:2009-Plateau}. A standard application of Allard's boundary regularity theorem exactly as in \cite[Chapter 4]{Eichmair:2009-Plateau}, using that the intersection $\partial M_- \cap \partial M_+$ is transverse, shows that the closure of the graph of $u_{t}^\epsilon$ in $M \times (- \frac{\epsilon}{t}, \frac{\epsilon}{t})$ is a $\C^{1, \alpha}$ manifold with boundary $(\partial M_- \cap \partial M_+) \times (-\frac{\epsilon}{t}, \frac{\epsilon}{t})$. The $\C^{1, \alpha}$ estimates near the boundary depend only on $C$ and the geometry of $\partial M_- \cap \partial M_+$; they are independent of $\epsilon, t >0$. 

We now pass the graphs of $u_{t}^\epsilon$ to a geometric subsequential limit as $\epsilon, t \searrow 0$. The existence and analysis of such limits is exactly as in \cite[Chapters 3, 4]{Eichmair:2009-Plateau} (which in turn are largely based on \cite{Schoen-Yau:1981-pmt2}). If in particular the limit along the subsequence $(t_n, \epsilon_n) \to (0, 0)$ were not a graph with the properties asserted in the theorem, there would be some $x \in \Omega$ such that the sequence $\{u_{t_n}^{\epsilon_n} (x)\}_{n=1}^\infty$ is unbounded. For definiteness, let us assume that $u_{t_n}^{\epsilon_n} (x) \to - \infty$, possibly after passing to a further subsequence. There exists a sequence of upward translations of the graphs of $u_{t_n}^{\epsilon_n}$ that converge, possibly after passing to a further subsequence, locally smoothly as hypersurfaces to a vertical cylinder $\Sigma \times \R$, where $\Sigma \subset \bar \Omega$ is a smooth properly embedded submanifold with boundary $\partial_- M \cap \partial_+ M$ that encloses a bounded region $\tilde \Omega$ with $\partial \Omega$ with  $x \in \tilde \Omega$, and such that $H_\Sigma + \tr_\Sigma(k) = 0$. Here, the mean curvature is computed with respect to the unit normal pointing out of $\tilde \Omega$. 
We can argue exactly as in \cite[Proposition 4.1 and Remark 4.1]{Eichmair:2010} that there exists a MOTS in $M_-$ that is homologous to $\partial M_-$ and which encloses $\{x\} \cup M_-$. (The point is that we can force a blow down of the Jang equation in the complement of the closure of $\tilde \Omega$ in $M_-$.) This contradicts the assumption that $\partial M_-$ is the outermost MOTS.   
\end{proof}
\end{theorem}


\section{Canonical blow up of the Jang equation} \label{sec:canonical}

In view of analogous results for Scherk-type minimal and constant mean curvature graphs on bounded domains, it is tempting to conjecture that the solutions of the Jang equation constructed in Proposition \ref{prop:existenceblowup} and Theorem \ref{thm:JangScherk} are unique with their properties. In Section \ref{sec:uniqueness} we prove such a uniqueness result in the special case where $k\equiv0$. In the case of general second fundamental form $k$, we will show in Theorem \ref{thm:maximalsolution} below that there exist \emph{canonical, pointwise maximal} solutions of the Jang equation for example when $M_- \subset M_+$. The proof of this result proceeds via a geometric variant of the Perron method that uses the outermost condition built into the definition of $M_+$ in lieu of a super solution for the problem. The basic ingredients are variations of classical PDE techniques, cf. in particular \cite{Jenkins-Serrin:1966, Jenkins-Serrin:1968, Serrin:1969, Serrin:1970} and also \cite{Eichmair:2009-Plateau} and the references therein, and the analysis of geometric limits of the Jang equation developed in \cite{Schoen-Yau:1981-pmt2}.

\begin{theorem} \label{thm:maximalsolution} Let $(M, g, k)$ be a complete asymptotically flat initial data set of dimension $n$, $3\leq n \leq 7$, fix an end, and let $M_+ \subset M$ be the complement of the total inner trapped region of $(M, g, k)$ with respect to that end. Assume that there exists a solution $u : M_+ \to \mathbb{R}$ of the Jang equation $\mc(u) + \tr(k)(u)=0$ such that $u(x) \to \infty$ as $\dist(x, \partial M) \to 0$ and such that $u(x) \to 0$ as $x \to \infty$ in the asymptotically flat end. The pointwise supremum of all such solutions is again a solution with the same properties. 
\begin{proof}
The results in Appendix \ref{sec:barriers} show that there exists $\Lambda \geq 1$ such that $|u(x)| \leq b_\Lambda (|x|)$ on $N= \{x \in M^+ : |x| > \Lambda\}$ for every solution $u$ of the Jang equation as in the statement of the theorem. 

Fix a smooth function $u : M_+ \to \R$ as in the statement of the theorem. Using translates of $u$ to obtain a priori oscillation bounds
 and standard methods as in \cite[Lemma 2.2]{Eichmair:2009-Plateau}, one sees that for every $x \in M_+$ there exists $\rho^D(x) \in (0, \dist(x, \partial M_+))$ small such that for every $\rho \in (0, \rho^D(x))$ the equation $\mc(v) + \tr(k)(v)=0$ on $B_\rho(x)$ with continuous boundary data on $S_\rho(x)$ admits a solution $v \in \C^{\infty} (B_\rho(x)) \cap \C^0(\bar B_\rho(x))$. The notion of (Perron) sub solutions $\underline u \in \C(M_+)$ and super solutions $\overline u \in \C(M_+)$ of the Jang equation $\mc(v) + \tr(k)(v) = 0$ on $M_+$ can thus be defined in the usual way. Consider the class of functions $\mathcal {S}_u = \{ \underline u \in \C (M_+) : \underline u \text{ is a Perron sub solution of the Jang equation, } \underline u \geq u \text{ on } M_+, \text{ and } |u(x)| \leq b_\Lambda (|x|) \text{ for all } x \in M_+ \text { with } |x| > \Lambda\}$. This class is closed under taking pointwise maximum and under lifting $\underline u \in \mathcal{S}_{u}$ to the function $\hat {\underline u} \in \C(M_+)$ that equals $\underline u$ on the complement of $B_\rho(x)$ and equals the solution $v$ of $H(v) + \tr(k)(v) = 0$ on $B_\rho(x)$ such that $v = u$ on $S_\rho(x)$, for every $r \in (0, \rho^D(x))$ and every $x \in M_+$. Note that $u \in \mathcal S_u$. The function $u^P : M_+ \to \R \cup \{\infty\}$ is defined pointwise by $u^P (x) = \sup_{\underline u \in \mathcal{S}_u} \underline u (x)$. Let $\Omega = \{x \in M_+ : u^P(x) < \infty\}$. Note that $N \subset \Omega$.

We claim that $\Omega$ is open, that $u^P|_\Omega$ is a smooth solution of the Jang equation, and that $\lim_{x \to y, x \in  \Omega} u^P(x) = \infty$ for every $y \in \partial \Omega$. To see this, fix $x \in \Omega$. Following the standard proof of the regularity of the Perron solution as in \cite[p. 25]{Gilbarg-Trudinger:1998}  we see that given $\rho \in (0, \rho^D(x))$ there exist $\{\underline u_i\}_{i=1}^\infty \subset \mathcal S_u$ such that 
$$u \leq \underline u_1 \leq \underline u_2 \leq \ldots \leq u^P$$
on  $B_\rho (x)$, such that 
$$H(\underline u_i) + \tr(k)(\underline u_i) = 0 $$ 
on  $B_\rho (x)$ for all $i = 1, 2, \ldots$, and such that 
$$\lim_{i \to \infty} \underline u_i (x) = u^P (x) < \infty.$$ 
The analysis of geometric limits of solutions of the Jang equation shows that the geometric limit of the graphs of $\underline u_i$ in $B_\frac{\rho}{2}(x) \times \R$ contains the graph of a smooth solution $\tilde u^{x, r} \in \C^\infty_{loc}(\Omega_0^{x, \rho})$ of the Jang equation above some open subset $\Omega_0^{x, \rho} \subset B_\rho(x)$. Moreover, we have that $\lim_{i \to \infty} \underline u_i(y) = \infty$ for all $y \in B_{\frac{\rho}{2}} (x) \setminus \Omega_0^{x, \rho}$, that $\lim_{z \to y, z \in \Omega_0^{x, \rho}} \tilde u^{x, \rho}(z) = \infty$ for all $y \in  B_{\frac{\rho}{2}}(x) \cap \partial \Omega_0^{x, \rho}$, and that $B_{\frac{\rho}{2}}(x) \cap \partial \Omega_0^{x, \rho}$ is a smooth properly embedded MITS in $B_{\frac{\rho}{2}}(x)$. (The point is that the functions $\underline u_i$ are bounded below by $u$ so that there can be no cylindrical components in their geometric limit, cf. the argument in Appendix \ref{sec:ige} and the properties listed in Step 4 in Subsection \ref{sec:OmegaNullEmpty}. Note that $\Omega_0^{x, \rho}$ might have several components.) Clearly, $\tilde u^{x, \rho} \leq u^P$ on $\Omega_0^{x, \rho}$. That $\tilde u^{x, \rho} = u^P$ on the connected component of $\Omega_0^{x, \rho}$ containing $x$ follows from the strong maximum principle for differences of solutions of the Jang equation, as in the standard proof of the regularity of Perron solutions. Since also $u^P \geq u$ everywhere on $M_+$, we can deduce all the properties of $\Omega$ and $u^P$ asserted at the beginning of this paragraph. 

The argument above shows that away from $\partial M_+$ the boundary of $\Omega$ is a smooth properly embedded MITS. That the boundary of $\Omega$ is smooth and embedded up to  $\partial_+ M$ follows from the characterization in Appendix \ref{sec:boundary} of boundaries of domains that support solutions of prescribed mean curvature equations with infinite boundary data. The definition of $M_+$ implies that $\Omega = M_+$. Clearly, $u^P$ has all the properties asserted in the conclusion of the theorem. 
\end{proof}
\end{theorem}

\begin{remark} By taking the smallest super solution instead of the largest sub solution in the proof of Theorem \ref{thm:maximalsolution} we obtain an analogous result where $M_+$ is replaced by $M_-$ and blow up to plus infinity at the boundary is replaced by blow down to minus infinity.  
\end{remark}


\section{Uniqueness of a blow up function $u$ at the outermost minimal surface} \label{sec:uniqueness}

Let $(M, g)$ be a complete asymptotically flat initial data set of dimension $n$, $2 \leq n \leq 7$, with $k \equiv 0$. Fix one of the ends and let $\Omega= M_- = M_+$ be as in the introduction. In this case, $\partial \Omega$ is called the horizon of $(M, g)$. Note that $\partial \Omega$ is a minimal surface. Let $\partial_- \Omega$ and $\partial_+ \Omega$ be unions of different components of $\partial \Omega$ such that $\partial \Omega = \partial_- \Omega \dot \cup \partial_+ \Omega$. \\

The arguments proving Proposition \ref{prop:existenceblowup} show that there exists a smooth solution $u : \Omega \to \R$ of the minimal surface equation $$\div \left( \frac{D u}{ \sqrt {1 + |D u|^2}} \right) = 0$$ such that $\lim_{x \to y, x \in \Omega} u(x) = \pm \infty$ for $y \in \partial_{\pm} \Omega$, and such that $u(x) \to 0$ as $|x| \to \infty$ in $\Omega$. Here we show that there is a solution with these properties. The proof is a straightforward adaption to our situation of a general argument due to J. Nitsche \cite{Nitsche:1965} as applied in e.g. \cite{Jenkins-Serrin:1966, Jenkins-Serrin:1968, Spruck:1972, Hauswirth-Rosenberg-Spruck:2009} to establish uniqueness of the Scherk-type graphs constructed there. We give the complete argument since the asymptotically flat ends require some care. \\

To see that $u$ is unique under the present assumptions, note first that the results in Appendix \ref{sec:boundary} show that the divergence of $u$ near $\partial \Omega$ is uniform in the distance to the respective components of the boundary, and that the upward and downward solutions of the graph converge geometrically to the vertical cylinders $\partial_+ \Omega \times \R$ and $\partial_- \Omega \times \R$ respectively. In particular, \begin{align} \label{eqn:convergencenormal} \lim_{x \to y, x \in \Omega} \frac{D u}{\sqrt {1 + |D u|^2}} (x) = \mp \nu(y)\end{align} for $y \in \partial_\pm \Omega$ where $\nu$ is the unit normal of $\partial \Omega$ pointing into $\Omega$. Using the argument in Appendix \ref{sec:barriers} we see that \begin{align} \label{eqn:decay} |u (x)| + |x| |D u(x)| = O(|x|^{2-\beta})\end{align} as $|x| \to \infty$ in $\Omega$ for every $\beta \in (2, n)$. \\

Suppose that $v : \Omega \to \R$ is a second solution with these properties. Fix $T \in (0, \infty)$ that is a regular value of both $(u-v)$ and $(v-u)$. Using (\ref{eqn:convergencenormal}) we see that
\begin{align*}
\lim _ {s \searrow 0} \int_{ \{x \in \Omega : \dist (x, \partial \Omega) = s, |(u - v)(x)| < T
\}} (u - v) g \left( \frac{D u}{\sqrt{1 + |D u|^2}} - \frac{D
v}{ \sqrt {1 + |D v|^2}}, \eta \right) d \mathcal{H}^{n-1}= 0. 
\end{align*}
Here, $\eta$ denotes the unit normal of $\{x \in \Omega : \dist (x, \partial \Omega) =
s\}$ pointing towards $\partial \Omega$. Using the decay estimates (\ref{eqn:decay}) for $u$ and $v$ we obtain that 
\begin{align*}
\lim_{r \to \infty} \int_{\{x \in \Omega : |x| = r, |u - v| < T\}}  (u - v) g \left(
\frac{D u}{\sqrt{1 + |D u|^2}} - \frac{D v}{\sqrt{1 + |D
v|^2}}, \eta\right)  d \mathcal{H}^{n-1}= 0
\end{align*}
where $\eta$ is the unit normal of $\{x \in \Omega: |x| = r\}$ pointing towards the
end. Using the divergence theorem and that $u, v$ satisfy the minimal surface equation, this implies that
\begin{equation*}
0 = \int_{\{x \in \Omega : |u - v| < T\}} g\left(D u - D v, \frac{D u
}{\sqrt{1 + |D u|^2}} - \frac{D v}{\sqrt{1 + |D v|^2}}\right) d
\mathcal{L}^n.
\end{equation*}
Using the strict convexity of the functions $\xi \to \sqrt{1 + g(\xi, \xi)}$ on $T_x M$ for all $x \in M$ we conclude that the integrand is pointwise non-negative with equality at $x \in \Omega$ if only if $D u = D v$ at $x$. 
It follows that $u$ and $v$ can only differ by a constant. Since we assume that they both tend to zero on the asymptotically
flat end, we obtain that $u = v$, as desired. \\

We see no way to extend this argument to non-zero $k$ at this point and have to contend with the existence of the canonical solution guaranteed under the hypotheses of Theorem \ref{thm:maximalsolution}. \\

\begin{remark} Let $u : \Omega \to \R$ be as above. It follows that $$\H^{n-1} (\partial_+ \Omega) - \H^{n-1} (\partial_-\Omega) = \lim_{r \to \infty} \int_{\{x \in \Omega : |x| = r\}} g (D u, D |x|) d \H^{n-1}.$$ Thus the unique blow up solution witnesses  the area of the horizon at infinity. \end{remark}


\section{Extensions of the classical Jenkins--Serrin theory} \label{sec:jss}

\subsection{Introduction}

The classical Jenkins--Serrin theory \cite{Jenkins-Serrin:1966, Jenkins-Serrin:1968} characterizes those bounded domains $\Omega \subset \R^2$ with piecewise smooth boundary for which there exists a solution $u : \Omega \to \R$ of the minimal surface equation such that 
\begin{align*}
\lim_{(x, y) \to (x_0, y_0) \in \partial_\pm\Omega} u(x, y) = \pm \infty
\end{align*}
and 
\begin{align*}
\lim_{(x, y) \to (x_0, y_0) \in \partial_0 \Omega} u(x, y) = \phi(x_0, y_0). 
\end{align*} 
Here, the sets $\partial_0 \Omega$, $\partial_- \Omega$, and $\partial_+ \Omega$ are unions of the smooth components of the boundary such that $\partial \Omega = \partial_+ \Omega \dot {\cup} \partial_- \Omega \dot {\cup} \partial_0 \Omega \dot \cup \{\text{corners}\}$, and $\phi \in C(\partial_0 \Omega)$ is a given function. A basic example of such a configuration is when $\Omega = (- \frac{\pi}{2}, \frac{\pi}{2}) \times (- \frac{\pi}{2}, \frac{\pi}{2})$ and $u(x, y) = \log \frac{\cos(x)}{\cos(y)}$. The graph of $u$ in $\R^3$ is of course the classical Scherk surface. In general, these domains are precisely those for which the connected components of $\partial_\pm \Omega$ are straight line segments such that no two segments in $\partial_- \Omega$ and no two segments in $\partial_+ \Omega$ have an endpoint in common, for which the geodesic curvature of $\partial_0 \Omega$ is non-negative, and which satisfy the Jenkins--Serrin condition \cite{Jenkins-Serrin:1966, Jenkins-Serrin:1968}. When $\partial_0 \Omega \neq \emptyset$ these conditions demand that the circumference of every polygon inscribed in $\Omega$ whose endpoints are chosen from the finitely many corner points is strictly greater than twice the total length of its sides that coincide with segments in $\partial_+ \Omega$ and also greater than twice the total length of its sides that coincide with segments in $\partial_- \Omega$. When $\partial_0 \Omega = \emptyset$ the Jenkins--Serrin condition is the same except for the inscribed polygon that is the whole domain; one demands that the length of $\partial_+ \Omega$ equals the length of $\partial_- \Omega$. \\

An important and influential development of the field was accomplished by J. Spruck \cite{Spruck:1972}, who has extended the classical Jenkins--Serrin theory to graphs $u : \Omega \subset \R^2 \to \R$ of constant mean curvature one. Provided that the piecewise smooth boundary $\partial \Omega$ consists of a union $\partial_+ \Omega$ of circular arcs of unit radius that are convex towards $\Omega$, a union of circular arcs $\partial_- \Omega$ of unit radius that are concave towards $\Omega$, and a union of boundary arcs $\partial_0 \Omega$ whose geodesic curvature is greater than or equal to one, he finds necessary and sufficient conditions for the existence of solutions $u$ assuming (arbitrary) continuous boundary values on $\partial_0 \Omega$ and tending to $\infty$ on approach towards $\partial_+ \Omega$ and $- \infty$ on approach towards $\partial_- \Omega$ \emph{provided} a sub solution $\underline u : \Omega^* \to \R$ \begin{align} \label{eqn:subsolution} \div \left( \frac{D \underline u}{\sqrt{1 + |D \underline {u}|^2}}\right) \geq 1\end{align} exists on the domain $\Omega^*$ obtained from flapping the negatively curved components of the boundary outward. (That this process gives a domain is a further additional assumption in \cite{Spruck:1972}.) The advantage of the piecewise smooth domain $\Omega^*$ over $\Omega$ is that its boundary arcs are all convex, so that solutions of the constant mean curvature equation with prescribed boundary values away from the corners can be constructed on it. \\

A further important contribution to the theory is due to U. Massari \cite{Massari:1977}, who has extended the Jenkins--Serrin theory to arbitrary dimensions and variable (Lipschitz) mean curvature in the case when $\partial_0 \Omega \neq \emptyset$. His techniques are different from J. Spruck's. In particular, no flapping of the boundary is required. In the special case where the mean curvature is constant, the necessary and sufficient conditions he provides for the existence of a solution are not obviously the same as those given in \cite{Spruck:1972}. We note that condition $(3.3)$ in \cite{Massari:1977} stands in, morally and technically, for the requirement (\ref{eqn:subsolution}) above that a sub solution of the prescribed mean curvature equation exists on the domain. In the book of E. Giusti \cite[Chapter 16]{Giusti:1984}, an extension of the technique of U. Massari to the ``slightly more complex" case where $\partial_0\Omega = \emptyset$ is presented in the minimal surface case. \\

More recently, the Jenkins--Serrin theory for minimal graphs has been extended to $\mathbb{H}^2 \times \R$ in \cite{Nelli-Rosenberg:2002} and to $M \times \R$ in \cite{Pinheiro:2009, Mazet-Rodriguez-Rosenberg:2011} (where $(M, g)$ is a general complete Riemannian surface). The Jenkins--Serrin--Spruck theory for constant mean curvature graphs has been developed for $\mathbb{H}^2 \times \mathbb{R}$ and $\mathbb{S}^2 \times \mathbb{R}$ in \cite{Hauswirth-Rosenberg-Spruck:2009} and in $\mathbb{M}^2 \times \R$ where $\mathbb{M}^2$ is a Hadamard surface by A. Folha and H. Rosenberg \cite{Folha-Rosenberg:2012}. A further important recent development are the results of P. Collin and H. Rosenberg \cite{Collin-Rosenberg:2010} and  A. Folha and S. Melo \cite{Folha-Melo:2011} who give necessary and sufficient conditions for the existence of Scherk-type minimal and constant mean curvature graphs on ideal polygons (with infinite area) in hyperbolic space. \\

As an example of a particularly interesting application of  Scherk-type graphs on Riemannian surfaces we mention the surprising construction of harmonic diffeomorphisms between the complex plane and the hyperbolic space by P. Collin and H. Rosenberg \cite{Collin-Rosenberg:2010}.  \\

In this section, we prove an extension of the Jenkins--Serrin--Spruck theory for domains $\Omega \subset M$ with $\partial_0 \Omega = \emptyset$  in Riemannian surfaces $(M, g)$ and for general $H_0 \in [0, \infty)$. When $H_0 = 0$, our result is marginally different from the corresponding result in \cite{Pinheiro:2009} in that we consider the possibility of closed geodesics in the blow up/blow down analysis (in Step $5$ below). In the case where $H_0 >0$, we neither hypothesize the existence of a sub solution on $\Omega$, nor do we make any additional assumptions regarding the existence of an auxiliary domain as in \cite{Spruck:1972, Hauswirth-Rosenberg-Spruck:2009}. \\

Our proof extends (for the most part verbatim) to appropriate domains in Riemannian manifolds of dimension $2 \leq n \leq 7$. The description of admissible domains in higher dimensions is cumbersome. Also, the Jenkins--Serrin--Spruck conditions become impossible to verify in examples that are not very symmetric. For this reason, we omit a detailed discussion of the extension to higher dimensions.  \\

In our proof we work with solutions of the (finite boundary value) Dirichlet problem that will in general not assume any particular boundary data but only have the right asymptotic behavior. This is an important difference with the construction in \cite{Spruck:1972} where solutions of the Dirichlet problem are constructed on an auxiliary domains $\Omega^*$ whose boundary is sufficiently convex away from finitely many points to construct solutions that take on prescribed continuous data. The construction of $\Omega^*$ in \cite{Spruck:1972} proceeds by flapping the negatively curved boundary components $\partial_- \Omega$ outward so they become convex arcs. This requires assumptions on the symmetry of the ambient manifold. For this reason, the construction in \cite{Spruck:1972, Hauswirth-Rosenberg-Spruck:2009} is carried out only in the simply connected space forms of dimension two. In \cite{Folha-Rosenberg:2012}, the existence of an extension of the original domain with properties similar to those of $\Omega^*$ is part of the assumption. See also \cite[Section 6]{Folha-Rosenberg:2012}. \\


\subsection {The case where $\partial_0 \Omega = \emptyset$} \label{sec:OmegaNullEmpty}

Let $(M, g)$ be a complete boundaryless Riemannian surface, let $H_0  \in [0, \infty)$, and let $\Omega \subsetneq M$ be a connected bounded open set such that $\partial \Omega = \partial \bar \Omega$. Here and below, we will use ``$\partial$" to denote the topological boundary of a set. We assume that $\partial \Omega$ is piecewise smooth and in fact the union of finitely many properly embedded arcs $\{A_i, B_j\}$ and properly embedded closed curves $\{E_k, F_l\}$ such that the outward geodesic curvature of each arc $A_i$ and closed curve $E_k$ is constant and equal to $H_0$ and such that the outward geodesic curvature of each arc $B_j$ and closed curve $F_l$ is constant and equal to $- H_0$. We assume that all the curves and the interior of all the arcs are pairwise disjoint, and that no two arcs $A_i$ and $A_{i'}$ and no two arcs $B_j$ and $B_{j'}$ have an endpoint in common. \\

The union of the closed curves $\{E_k\}$ and the interiors of the positively curved arcs $\{A_i\}$  is denoted by $\partial_+ \Omega$. The union of the closed curves $\{F_l\}$ and the interiors of the negatively curved arcs $\{B_j\}$  is denoted by $\partial_-\Omega$.  The endpoints of the arcs $\{A_i, B_j\}$ are called the \emph{corners} of $\partial \Omega$. The assumptions imply (via the strong maximum principle when $H_0 >0$) that any two arcs that share an endpoint meet at a non-zero angle. \\

A \emph{generalized polygon} is a non-empty open subset $P\subset \Omega$ with $\partial P = \partial \bar P$ and such that $\partial P$ is piecewise smooth and consists of finitely many of the following building blocks: 
\begin{enumerate} [(i)]
\item Finitely many arcs of constant geodesic curvature $
H_0$ whose endpoints are amongst the corners of $\partial\Omega$ and whose interiors are embedded and pairwise disjoint. We also require that each arc whose interior intersects $\partial \Omega$ is one of $\{A_i, B_j\}$. 
\item Pairwise disjoint embedded closed curves of constant geodesic curvature $\pm H_0$ that either lie entirely in $\Omega$ or coincide with one of $\{E_k, F_l\}$ and which are disjoint from the boundary arcs.
\end{enumerate}

\begin{theorem} \label{thm:mainjss} Let $(M, g)$ and $\Omega \subset M$ be as above. A necessary and sufficient condition for the existence of a smooth function $u: \Omega \to \R$ such that 
\begin{align} \label{eqn:jssgraph}
\div \left (\frac{D u}{ \sqrt{1 + |D u|^2}} \right) = H_0
\end{align}
with 
\begin{align*} \lim_{x \to x_0, x \in \Omega} u(x) =  \left \{ \begin{array}{ll}  \infty & \text{if } x_0 \in \partial_+ \Omega \\ - \infty & \text{if } x_0 \in \partial_- \Omega \end{array} \right. 
\end{align*}
is that \begin{align} \label{eqn:totalflux}  \H^1_g(\partial_+ \Omega) =  H_0 \L^2_g(\Omega) + \H^1_g(\partial_- \Omega) \end{align} and that  \begin{align} \label{eqn:partialfluxA} 2 \H^1_g(\partial_+ \Omega \cap \partial P) < \H^1_g(\partial P) + H_0 \L^2_g(P)
\end{align}
and 
\begin{align} \label{eqn:partialfluxB} 2 \H^1_g(\partial_- \Omega \cap \partial P) < \H^1_g (\partial P) - H_0 \L^2_g(P)\end{align} for every generalized polygon $P \subsetneq \Omega$. 
\end{theorem}

The conditions (\ref{eqn:totalflux}), (\ref{eqn:partialfluxA}), (\ref{eqn:partialfluxB}) appear in the classical work of Jenkins--Serrin (when $H_0 = 0$ and $(M, g)$ is $\R^2$ with the Euclidean metric) and its generalization due to J. Spruck (when $H_0 >0$ and $(M, g)$ is Euclidean space), and also in the work of L. Hauswirth, H. Rosenberg, and J. Spruck \cite{Hauswirth-Rosenberg-Spruck:2009} (where $H_0 > 0$ and $(M, g)$ is one of $\R^2, \mathbb{S}^2, \mathbb{H}^2$ with their constant curvature metrics). The necessity of these conditions follows from a standard argument, that we summarize briefly: \\

Let $u : \Omega \to \R$ be as in the statement of Theorem \ref{thm:mainjss}. Let $U \subset M$ be a non-empty and open subset  such that $U \cap \partial_\pm \Omega = U \cap \partial \Omega$. The discussion in Appendix \ref{sec:boundary} shows that the graphs of the functions $u \mp t$ converge as hypersurfaces smoothly on compact subsets of $U \times \R$ to $\partial_\pm \Omega \times \R$ as $t \to \infty$. In particular, as $x \in \Omega$ approaches a point $x_0 \in \partial_\pm \Omega$, the horizontal part of the downward unit normal of these graphs, $X=  (1 + |D u|^2)^{-1/2} Du$, converges to $\pm$ the outward pointing unit normal of $\partial_\pm \Omega$ at $x_0$.  The necessity of condition (\ref{eqn:totalflux}) follows from applying the divergence theorem to the vector field $X$ on smooth interior approximations of the domain $\Omega$. The necessity of conditions (\ref{eqn:partialfluxA}), (\ref{eqn:partialfluxB}) follows from the same argument applied to generalized polygons $P \subsetneq \Omega$, using also that $|X (x)| < 1$ for $x \in \partial P \cap \Omega$.  \\

\begin{remark} A consequence of the existence of a graph as in (\ref{eqn:jssgraph}) is that for every $\gamma \in \{A_i, B_j, E_k, F_l\}$ we have that 
$$\int_\gamma |\bar D \psi|^2  \geq \int_\gamma (H_0^2 + \kappa_g) \psi^2$$ 
for every $\psi \in C^1(\gamma)$ that vanishes near the boundary of $\gamma$. Here, $\bar D$ is the (tangential) gradient and $\kappa_g$ is half the scalar curvature of $(M, g)$. This implies, for example, that in Euclidean space and when $H_0 = 1$, a domain that satisfies the Jenkins--Serrin--Spruck conditions (\ref{eqn:totalflux}), (\ref{eqn:partialfluxA}), (\ref{eqn:partialfluxB}) must also satisfy $\H^1_\delta (A_i), \H^1_\delta(B_j) \leq \pi$. The condition that $\H^1_\delta (B_j) < \pi$ was part of the assumptions in \cite[p. 16]{Spruck:1972}. In fact, our proof (see property (\ref{item10}) in Step 4 below) shows that the conditions (\ref{eqn:totalflux}), (\ref{eqn:partialfluxA}), (\ref{eqn:partialfluxB}) need only be verified for generalized polygons whose boundary components are stable. This sharpening of the classical Jenkins--Serrin--Spruck condition is useful when constructing examples.
\end{remark}

\begin{remark}  The complement of a generalized polygon $P$ in $\Omega$ is again a generalized polygon. Condition (\ref{eqn:partialfluxB}) for a generalized polygon $P \subsetneq \Omega$ follows from condition (\ref{eqn:partialfluxA}) applied to $\Omega \setminus \bar P$ in view of (\ref{eqn:totalflux}).
\end{remark}


\noindent {\bf Step 1: Construction of the auxiliary domain $\hat \Omega$}

Fix a component $\gamma$ of $\partial_\pm \Omega$ and a smooth function $\Theta_\gamma \in \C^\infty(\bar \gamma)$ that is positive on $\gamma$ and which vanishes on its (manifold) boundary. Let $\nu$ be the unit normal of $\gamma$ pointing out of $\Omega$. The piecewise smooth domain $\hat \Omega$ is obtained from $\bar \Omega \setminus \{\text{vertices}\}$ by adding the \emph{crescents}  
$$\text{Cr}_\gamma = \{ \exp_\theta (t  \Theta_\gamma(\theta) \nu (\theta) ) : t \in (0, \epsilon) \text{ and } \theta \in \text{int}(\gamma) \}$$
as $\gamma$ ranges over all components of $\partial_\pm \Omega$. Here, $\epsilon > 0$ small is chosen so that there are no issues with the regularity of the exponential map and such that the crescents $\text{Cr}_\gamma$  are pairwise disjoint.  \\


\noindent {\bf Step 2: Construction of barriers on $\text{Cr}_\gamma$ and the functions $H_k(x)$} 

Let $\gamma$ be a component of $\partial_+ \Omega$. We would like to find solutions that tend to $\infty$ on approach to $\gamma$ and hence require a sub solution.  The hypersurface 
$$\{ (\exp_\theta (\epsilon  e^h \Theta_\gamma (\theta) \nu (\theta)), h) : h \in (- \infty, 0) \text{ and } \theta \in \text{int} (\gamma)\}$$ 
of $M \times \R$ is the graph of a (locally) smooth function $\underline u_\gamma : \text{Cr}_\gamma \to \R$ whose downward unit normal corresponds to the outward unit normal of the cylinder. \\

Let $\gamma$ be a component of $\partial_- \Omega$. We would like to find solutions that tend to $- \infty$ (and hence require a super solution). As above, the hypersurface 
$$\{ (\exp_\theta (\epsilon e^{-h} \Theta_\gamma (\theta) \nu (\theta)), h) : h \in (0, \infty) \text{ and } \theta \in \text{int} (\gamma) \}$$ 
is the vertical graph of a locally smooth function $\overline u_\gamma : \text{Cr}_\gamma \to \R$.\\

The function $H : \hat \Omega \to \R$ defined as 
\begin{align*}
H(x) =  \left \{ \begin{array}{lll}  H_0  & \text{if } x \in \bar \Omega \setminus \{ \text{corners}\},      \\ \div \left( \frac{D \overline u_\gamma }{\sqrt{1 + |D \overline u_\gamma|^2}}\right) (x)  & \text{if } x \in \text{Cr}_\gamma \text{ and if } \gamma \text{ is a component of } \partial_- \Omega, \text{ and}\\ \div \left( \frac{D \underline u_\gamma }{\sqrt{1 + |D \underline u_\gamma|^2}}\right) (x)  & \text{if } x \in \text{Cr}_\gamma \text{ and if } \gamma \text{ is a component of } \partial_+ \Omega  \end{array} \right. 
\end{align*}
is locally Lipschitz. Moreover, $H(x) = H_0 + O(\epsilon)$ uniformly on $\hat \Omega$. \\

Let $\chi : \hat \Omega \to [-1, 1]$ be a locally smooth function such that $\chi \equiv \pm 1$ near $\text{Cr}_\gamma$ for $\gamma \in \partial_\pm \Omega$. Let $k \geq 1$. We define a locally Lipschitz function $H_k$ on $\hat \Omega$ by $H_k (x) = H (x) -  k^{-1/2} \chi (x)$. Note that $\sqrt{k} + \underline u_\gamma$ is a sub solution of the equation 
\begin{align} \label{eqn:regularized} 
\div \left( \frac{D u}{\sqrt{1 + |D u|^2}} \right)=  H_k + \frac{1}{k} \, u
\end{align} 
on $\text{Cr}_\gamma$ when $\gamma$ is a component of $\partial_+ \Omega$, and that $- \sqrt{k} + \overline u_\gamma$ is a super solution for this equation on $\text{Cr}_\gamma$ when $\gamma$ is a component of $\partial_- \Omega$. \\

Fix a constant $C > \sup_{k \geq 1, x \in \hat \Omega} |H_k(x)|$. Then $- C k$ and $Ck$ are, respectively, sub and super solutions for (\ref{eqn:regularized}) on $\hat \Omega$. \\

The introduction of the capillarity regularization (\ref{eqn:regularized}) of the prescribed mean curvature equation so that  large constants become barriers is exactly as in \cite{Schoen-Yau:1981-pmt2}. \\


\noindent {\bf Step 3: The construction of $u_k$}

Let $u_k \in \C^{2, \alpha}_{loc} (\hat \Omega)$ be the largest (Perron) sub solution of equation (\ref{eqn:regularized})
that lies below $C k$ on all of $\hat \Omega$ and below $- \sqrt{k} + \overline u_\gamma$ on all crescents $\text{Cr}_\gamma$ corresponding to components $\gamma$ of $\partial_- \Omega$. To justify the existence of such a solution, we refer to \cite[p. 375]{Serrin:1970}, the interior gradient estimate stated in Appendix \ref{sec:ige}, and also \cite[Theorems 1.1 and 1.4]{Spruck:2007}, \cite[Chapter 16]{Gilbarg-Trudinger:1998}, and \cite[Chapter 3]{Eichmair:2009-Plateau}. The maximum principle implies that $u_k$ lies above  $- C k$ on all of $\hat \Omega$ and above $\sqrt{k} + \underline{u}_\gamma$ on all crescents $\text{Cr}_\gamma$ corresponding to components $\gamma$ of $\partial_+ \Omega$. From (\ref{eqn:regularized}) we see that the mean curvature of the graph of $u_k$ is bounded uniformly by $2C$ on $\hat \Omega$. \\


\noindent {\bf Step 4: Geometric limits of $\graph(u_k)$}

We claim that there is a subsequence $\{u_{k_i}\}$ of $\{u_k\}$ and there exist disjoint open subsets $\Omega_0, \Omega_+, \Omega_-$ of $\hat \Omega$ and $u \in C^{2, \alpha}_{loc} (\Omega_0)$ with the following properties: 

\begin{enumerate} [(a)]
\item \label{item0} $\hat \Omega = (\overline \Omega_0 \cup \overline {\Omega}_- \cup \overline{\Omega}_+) \cap \hat \Omega$. In particular, the topological boundaries $\partial \Omega_0$, $\partial \Omega_-$, $\partial \Omega_+$ of $\Omega_0$, $\Omega_-$, $\Omega_+$ locally separate $\Omega_0$, $\Omega_-$, $\Omega_+$ from their respective complements $\hat \Omega \setminus \Omega_0$,  $\hat \Omega \setminus \Omega_-$, $\hat \Omega \setminus \Omega_+$ in $\hat \Omega$. Moreover, $\partial \Omega_0 \cap \hat \Omega$, $\partial \Omega_- \cap \hat \Omega$,  and $\partial \Omega_+ \cap \hat \Omega$ are properly embedded $\C^{2, \alpha}$ hypersurfaces in $\hat \Omega$. 
\item \label{item1} For every $x \in \Omega_+$ there exists an open neighborhood of $x$ in $\Omega_+$ so that $u_{k_i} (y)$ exceeds a given constant for all $y$ in this neighborhood, provided $i$ is sufficiently large.  Put differently, $u_{k_i}$ diverges to plus infinity locally uniformly on $\Omega_+$. 
\item For every $x \in \Omega_-$ there exists an open neighborhood of $x$ in $\Omega_-$ so that $u_{k_i} (y)$ lies below a given constant for all points $y$ in this neighborhood, provided $i$ is sufficiently large. Put differently, $u_{k_i}$ diverges to minus infinity locally uniformly on $\Omega_-$. 
\item We have that $u_{k_i} \to u$ in $C^{2, \alpha}_{loc}(\Omega_0)$. In particular, $$\div \left(\frac{D u}{ \sqrt{1 + |D u|^2}}\right) =  H$$ on $\Omega_0$.  
\item \label{item2} $\text{Cr}_\gamma \subset \Omega_\pm$ when $\gamma$ is a component of $\partial_\pm \Omega$. 
\item The sets $\partial \Omega_0 \cap (\hat \Omega \cup \{\text{corners}\})$, $\partial \Omega_- \cap (\hat \Omega \cup  \{\text{corners}\})$, and $\partial \Omega_- \cap (\hat \Omega  \cup \{\text{corners}\})$ consist of finitely many arcs and closed curves of constant geodesic curvature in $\overline {\Omega}$. These arcs and closed curves are pairwise disjoint in $\Omega$. The arcs are properly immersed and embedded in $\Omega$. Their endpoints are corners of $\Omega$, and the endpoints of any one arc may coincide. The closed curves are contained in $\Omega$ and they are properly embedded. 
\item \label{item11} If $\gamma$ is a component of $\partial \Omega_\pm \cap \partial \Omega_0$, then $\lim_{x \in \Omega_0, x \to x_0} u (x) = \pm \infty$ uniformly near $x_0 \in \text{int} (\gamma)$. 
\item \label{item5} The geodesic curvature of a component $\gamma$ of $\partial \Omega_0 \cap \partial \Omega_+$ is constant and equal to $H_0$ when we orient $\gamma$ by the unit normal $\nu$  pointing into $\Omega_+$. Every divergent series of downward translations of the hypersurface $\text{graph} (u) = \{(x, u(x)) : x \in \Omega_0\}$ converges to $(\partial \Omega_0 \cap \partial \Omega_+)\times \R$ in $\C^{2, \alpha}$ on compact subsets of $\hat \Omega \times \R$. 
\item  \label{item6} The geodesic curvature of a component $\gamma$ of $\partial \Omega_0 \cap \partial \Omega_-$ is constant and equal to $H_0$ when we orient $\gamma$ by the unit normal $\nu$  pointing into $\Omega_0$. Every divergent series of upward translations of the hypersurface $\text{graph} (u) = \{(x, u(x)) : x \in \Omega_0\}$ converges to $(\partial \Omega_0 \cap \partial \Omega_-)\times \R$ in $\C^{2, \alpha}$ on compact subsets of $\hat \Omega \times \R$. 
\item \label{item7} The geodesic curvature of a component $\gamma$ of $\partial \Omega_- \cap \partial \Omega_+$ is constant and equal to $H_0$ when we orient $\gamma$ by the unit normal $\nu$  pointing into $\Omega_+$. There exists a cylindrical neighborhood of $\gamma \times \R$ in $\hat \Omega \times \R$ in which the graphs $\{(x, u_{k_i} (x)) : x \in \hat \Omega\}$ converge in $\C^{2, \alpha}$ to $\gamma \times \R$ on compact subsets.
\item \label{item16} The graphs $\{(x, u_{k_i} (x)) : x \in \hat \Omega\}$ converge as embedded $\C^{2, \alpha}$ hypersurfaces on compact subsets of $\hat \Omega \times \R$ to the union of the cylinders $(\Omega \cap \partial_0 \Omega)\times \R$, $(\partial_- \Omega \cap \partial_+ \Omega) \times \R$ and the graph $\{(x, u(x)) : x \in \Omega_0\}$, as $i \to \infty$. 
\item \label{item9} The vector fields $D u_k /\sqrt {1 + |D u_k|^2}$ are locally equicontinuous in $\hat \Omega$. 
\item \label{item8} With $\gamma$ and $\nu$ as in (\ref{item5}) - (\ref{item7}) we have that $$\lim_{i \to \infty} \int_\gamma \frac {g(\nu, D u_{k_i})} {\sqrt{1 + |D u_{k_i}|^2}} d \H^1_g = \H^1_g (\gamma).$$ In fact, the integrand on the left converges to $1$ locally uniformly on $\text{int}(\gamma)$.
\item \label{item10} The arcs and closed curves $\gamma$ in (\ref{item5}), (\ref{item6}), and those in (\ref{item7}) interior to $\Omega$ are stable in the sense that 
\begin{align*} 
 \int_\gamma (H^2_0 + \kappa _g) \psi^2 d \H^1_g \leq \int_\gamma |\bar D \psi|^2 d \H^1_g
 \end{align*} for all $\psi \in \C^1(\gamma)$ with $\text{supp} (\phi) \subset \text{int}(\gamma)$.
Here, $\kappa_g$ is half the scalar curvature of $(M, g)$ and $\bar D$ is the (tangential) gradient of $\psi$. 
\end{enumerate}
The properties listed here extend classical results about limits of monotone sequences in the Jenkins--Serrin--Spruck theory, cf. \cite{Jenkins-Serrin:1966, Jenkins-Serrin:1968, Spruck:1972, Hauswirth-Rosenberg-Spruck:2009, Pinheiro:2009}. For limits of not necessarily monotone sequences, some of these properties can be inferred directly from the results of L. Mazet \cite{Mazet:2007}. The ideas in \cite{Mazet:2007} have been employed to prove Jenkins-Serrin-Spruck type results for minimal graphs supported on domains in Riemannian surfaces in \cite{Mazet-Rodriguez-Rosenberg:2011} and for constant mean curvature graphs in \cite{Folha-Melo:2011}. 

The above properties are also variations of classical results on generalized solutions of the minimal surface equation \cite{Giusti:1978, Massari-Miranda:1984} or the geometric theory of Jang equation \cite{Schoen-Yau:1981-pmt2}. We provide a few details below to assist the reader. \\

The barriers used in the construction of $u_k$ in Step 4 can be used to justify (\ref{item2}). \\

Properties (\ref{item0}) - (\ref{item8}) can be deduced from the compactness and regularity properties of almost minimizing boundaries along with the geometric Harnack principle (as in Appendix \ref{sec:ige}), which ensures that geometric limits of our graphs are made up of graphical and cylindrical components. For limits of minimal graphs, this is the approach of \cite{Massari-Miranda:1984}. For the case of geometric limits of solutions of the regularized Jang equation (including the capillarity regularization), this has been worked out in detail in \cite{Eichmair:2009-Plateau}. We refer the reader to \cite{Eichmair:2009-Plateau} for statements of more general results and references. \\

Property (\ref{item9}) follows from Lemma \ref{lem:equicontinuity} in Appendix \ref{sec:equicontinuous}. \\

The argument leading to (\ref{item10}) is very similar to that in Appendix \ref{sec:boundary}. The Jacobi identity (\ref{Jacobiidentity}) is replaced by the differential inequality 
\begin{align*} \Delta_{G_k} \nu^{3}_k + (|h_k|^2 + \Ric_{g + dx^3 \otimes dx^3}(\nu_k, \nu_k)) \nu^{3}_k &= - g \left( \frac{ D u_k } {{\sqrt{1 + |D u_k|^2}}},  D (H_k + \frac{1}{k} \, u_k)\right) \nu^3_k \\ & \leq k^{-1/2} |D \chi| \nu^3_k
\end{align*}
valid on $\Omega$ where all geometric quantities on the left are computed for the graph $G_k = \{(x, u_k(x)) : x \in \Omega\}$. The additional contribution to the stability inequality (\ref{eqn:stabilityG}) disappears when we take geometric subsequential limits of $G_k$ and its vertical translates as $n \to \infty$. Since the arcs $\gamma$ for which (\ref{item10}) is asserted appear as cross-sections of vertical cylinders that appear in such limits, and because $|h|^2 + \Ric_{g+dx^3 \otimes dx^3} (\nu, \nu)$ reduces to $H_0^2 + \kappa_g$ on such cross-sections, we are done. \\


\noindent {\bf Step 4: Analysis of the limit using the Jenkins--Serrin--Spruck conditions}

The analysis of the geometric limit of the graphs of the solutions $u_k : \Omega \to \R$ using the Jenkins--Serrin--Spruck condition below is in many ways similar to that in \cite{Spruck:1972} (see in particular Sections 5 and 6 therein) or \cite[Section 7]{Hauswirth-Rosenberg-Spruck:2009}. However, because we do not assume the existence of a sub solution for the original equation, our main technical step, Case b below, is quite different.  \\

The properties listed in Step 4 show that the components $P$ of $\Omega_0 \cap \Omega$,  $\Omega_+\cap \Omega$, and $\Omega_-\cap \Omega$ are generalized polygons in $\Omega$. If $P$ is a component of $\Omega_- \cap \Omega$, then 
\begin{align} \label{eqn:flux-}
 H_0 \L^2_g(P) +  \limsup_{i \to \infty} \int_P \frac{u_{k_i}}{k_i} d \L^2_g \geq \\ \H^1_g (\partial P \cap \partial_+ \Omega) +  \H^1_g (\partial P \cap \Omega)  +  \liminf_{i \to \infty} \int_{\partial P \cap \partial_- \Omega} \frac {g(\nu, D u_{k_i})} {\sqrt{1 + |D u_{k_i}|^2}}d \H^1_g. \nonumber
\end{align}
Here, $\nu$ is the unit normal pointing out of $\Omega$. The second term on the left is always non-positive. Similarly, if $P$ is a component of $\Omega_+ \cap \Omega$, then 
\begin{align} \label{eqn:flux+}
 H_0 \L^2_g(P)  +  \liminf_{i \to \infty} \int_P \frac{u_{k_i}}{k_i}  d \L^2_g \leq  \\  - \H^1_g (\partial P \cap \partial_- \Omega) -  \H^1_g (\partial P \cap \Omega) + \limsup_{i \to \infty} \int_{\partial P \cap \partial_+ \Omega}  \frac {g(\nu, D u_{k_i})} {\sqrt{1 + |D u_{k_i}|^2}}d \H^1_g, \nonumber
\end{align}
where again the unit normal $\nu$ points out of $\Omega$. The second term on the left is always non-negative. \\

We consider the following cases: \\

\noindent {\bf Case a: $\emptyset \neq \Omega_- \cap \Omega \subsetneq \Omega$}

Let $P$ be a component of $\Omega_- \cap \Omega$. The Jenkins--Serrin--Spruck condition for $P$ implies that 
\begin{align} \nonumber
 H_0 \L^2_g(P) < \H^1_g (\partial P \cap \partial_+ \Omega) +  \H^1_g (\partial P \cap \Omega) -  \H^1_g (\partial P \cap \partial_- \Omega). 
\end{align} 
This contradicts our assumption (\ref{eqn:flux-}), since clearly $$\left| \int_{\partial P \cap \partial_- \Omega}  \frac {g(\nu, D u_{k_i})} {\sqrt{1 + |D u_{k_i}|^2}} d \H^1_g \right| \leq \H^1_g(\partial P \cap \partial_- \Omega).$$ Thus Case a cannot occur. \\

\noindent {\bf Case a': $\emptyset \neq \Omega_+ \cap \Omega \subsetneq \Omega$} 

Let $P$ be a component of $\Omega_+ \cap \Omega$. The Jenkins--Serrin--Spruck condition for $P$ implies that 
\begin{align} \nonumber
H_0 \L^2_g(P) > \H^1_g (\partial P \cap \partial_+ \Omega) -  \H^1_g (\partial P \cap \Omega) -  \H^1_g (\partial P \cap \partial_- \Omega), 
\end{align} 
contradicting (\ref{eqn:flux+}), since clearly $$\left| \int_{\partial P \cap \partial_+ \Omega}  \frac {g(\nu, D u_{k_i})} {\sqrt{1 + |D u_{k_i}|^2}} d \H^1_g \right| \leq \H^1_g(\partial P \cap \partial_+ \Omega).$$ Thus Case a' cannot occur. \\

\noindent {\bf Case b: $\Omega \subset \Omega_-$}

The Jenkins--Serrin--Spruck condition for $P = \Omega_- \cap \Omega = \Omega$ implies that 
\begin{align} \nonumber
 H_0 \L^2_g(\Omega) = \H^1_g (\partial_+ \Omega) -  \H^1_g (\partial_- \Omega). 
\end{align}
In conjunction with (\ref{eqn:flux-}) we conclude that  $$\limsup_{i \to \infty} \int_P \frac{u_{k_i}}{ {k_i}}d \L^2_g = 0$$ and that $$\liminf_{i \to \infty} \int_{\partial_- \Omega}  \frac {g(\nu, D u_{k_i})} {\sqrt{1 + |D u_{k_i}|^2}} d \H^1_g = - \H^1_g (\partial_- \Omega).$$ Passing to a further subsequence, if necessary, we see that 
\begin{align} \label{eqn:conclusion1Caseb}
\lim_{i \to \infty} \L^2_g (\{ x \in \Omega : \frac{u_{k_i}}{{k_i}} < - \epsilon\}) = 0
\end{align}
for every $\epsilon > 0$
and, using (\ref{item9}), that 
\begin{align} \label {eqn:conclusion2Caseb}
\lim_{i \to \infty}  \frac {g(\nu, D u_{k_i})} {\sqrt{1 + |D u_{k_i}|^2}} = -1
\end{align} locally uniformly on $\partial_- \Omega$
where $\nu$ is the unit normal pointing out of $\Omega$. \\

Fix a component $\gamma$ of $\partial_-\Omega$ and let $z \in \gamma$. It follows from the assumptions that $z \in \Omega_-$. Consider the functions $\tilde u_{k_i} (x) = u_{k_i} (x) - u_{k_i} (z)$. (This is an upward translation for $k_i$ large.) Then
$$\div \left( \frac{D \tilde u_{k_i}}{\sqrt{1 + |D \tilde u_{k_i}|^2}} \right)=  H_{k_i} + \frac{\tilde u_{k_i}}{k_i} + \frac{u_{k_i}(z)}{k_i}.$$ Using that $ |u_{k_i} (x) | \leq C k_i$ for all $x \in \hat \Omega$ we see that the mean curvature of these graphs is uniformly bounded. We pass to a further subsequence so that $k_i^{-1}  u_{k_i} (z)$ converges to a constant $c \in [- C, 0]$. We pass to a further subsequence so that the graphs of the $\tilde u_{k_i}$ converge geometrically in $\C^{2, \alpha}$ to a union of properly embedded graphs and cylinders on compact subsets in $\hat \Omega \times \R$. The mean curvature of these graphs and cylinders at a point $(x, x^{3}) \in \hat \Omega \times \R$ in the geometric limit is $H(x) + c$. The point $(z, 0) \in \gamma \times \R$ is contained in the geometric limit. Using (\ref{eqn:conclusion2Caseb}) we see that $(z, 0)$ is contained in a cylindrical component $\tilde \gamma \times \R$ of the limit, where $\tilde \gamma \subset \hat \Omega$ is a properly embedded curve whose mean curvature at $x \in \gamma$ is given by $H(x) + c$. Moreover, the tangent spaces of $\gamma$ and $\tilde \gamma$ agree together with their orientation at any point of $\gamma \cap \tilde \gamma$. \\

We claim that $\gamma = \tilde \gamma$. In particular, $c = 0$. To see this, we distinguish two cases. \\

First, assume that $\tilde \gamma \subset \overline {\text{Cr}}_\gamma$. Then the assertion is a consequence of the maximum principle. (We use that $c \leq 0$ here.) \\

Second, assume that $\tilde \gamma \cap \Omega \neq \emptyset$. Recall that every $y \in \gamma$ has an open neighborhood in $\hat \Omega$ that is separated by $\gamma$ into two components such that $\tilde u_{k_i}$ tends to plus infinity locally uniformly in one component as $k \to \infty$, and such that $\tilde u_{k_i}$ tends to minus infinity locally uniformly in the other component. Since $\tilde \gamma \cap \Omega \neq \emptyset$, we conclude that 
$$\{x \in \Omega : \limsup_{i \to \infty} k_i^{-1 } \tilde u_{k_i}(x) \leq 0\} = \{x \in \Omega : \limsup_{i \to \infty}  k_i^{-1} u_{k_i}(x) \leq c\}$$ 
contains a non-empty open subset. In conjunction with (\ref{eqn:conclusion1Caseb}) we conclude that $c=0$. It follows that in this case, $\gamma$ and $\tilde \gamma$ satisfy the same geometric equation. Further, we know that they intersect non-trivially, and that at any point of intersection they intersect tangentially with the same orientation. The Hopf boundary point lemma shows that $\gamma = \tilde \gamma$, which contradicts the assumption that $\tilde \gamma \cap \Omega \neq \emptyset$. It follows that $\gamma = \tilde \gamma$.\\

The argument in the preceding paragraph shows that there exists a relatively open neighborhood $U_\gamma$ of $\gamma$ in $\hat \Omega$ that is disjoint from the crescents corresponding to the positively curved boundary components such that $\tilde u_{k_i}$ converges to $- \infty$ locally uniformly in $U_\gamma \cap \text{Cr}_\gamma$, and to $\infty$ in $U_\gamma \cap \Omega$. \\

We can repeat the above reasoning for any of the components $\gamma_1, \ldots, \gamma_m$ of $\partial_- \Omega$, choosing a point $z_i \in \gamma_i$ for each component. Passing to a further subsequence and relabeling if necessary, we may assume that $u_{k_i} (z) \geq u_{k_i} (z_i)$ for all $i \in \{1, \ldots, m\}$ where $z = z_1$. Also, we may pick $y \in \Omega$ near $\gamma_1$ so that 
$$0 \geq \frac{u_{k_i} (y)}{k_i} \geq \frac{u_{k_i} (z)}{k_i}  \to 0.$$ 
Finally, let $\hat u_{k_i} = u_{k_i} - u_{k_i} (y)$. For $k$ large, this is an upward translation. Then   
$$\div \left( \frac{D \hat u_{k_i}}{\sqrt{1 + |D \hat u_{k_i}|^2}} \right)=  H_{k_i} + \frac{\hat u_{k_i}}{k_i} + \frac{u_{k_i}(y)}{k_i}.$$
For the remainder of the argument, we replace $\text{Cr}_\gamma$ by $U_\gamma \cap \text{Cr}_\gamma$ (keeping the same notation) for all components $\gamma$ of $\partial_- \Omega$. This may shrink the auxiliary domain $\hat \Omega$ slightly. We have that $\tilde u_{k_i} (x) \to - \infty$ locally uniformly in these new crescents $\text{Cr}_\gamma$. Note that $\tilde u_{k_i} \to \infty$ locally uniformly in $\text{Cr}_\gamma$ when $\gamma$ is a component of $\partial_+\Omega$. We can take a subsequential geometric limit of the graphs of $\tilde u_{k_i}$ just as we did for the original sequence $u_{k_i}$ so that (\ref{item0}) - (\ref{item9}) continue to hold. We use $\hat \Omega_0, \hat \Omega_-, \hat \Omega_+$ instead of $\Omega_0, \Omega_-, \Omega_+$ to avoid confusion. As before, the components of $\Omega \cap \hat \Omega_0$ and $\Omega \cap \hat \Omega_\pm$ are generalized polygons. We claim that $\Omega \subset \hat \Omega_0$. To see this, note that $(y, 0)$ is contained in the geometric limit. This implies that $y \in \partial \hat \Omega_- \cup \partial \hat \Omega_+ \cup \hat \Omega_0$ so that $\Omega \subset \hat \Omega_\pm$ is impossible. The cases $\emptyset \neq \Omega \cap \hat \Omega_\pm \subsetneq \Omega$ can be ruled out exactly as in Cases a and a' above. It follows that $\Omega = \hat \Omega_0$, and we can conclude as in Case c below. \\

\noindent {\bf Case b': $\Omega \subset \Omega_+$}

Exactly as in Case b, we can conclude that a sequence of downward translations of $u_{k_i}$ will converge to a solution of the original problem. We point out that the analysis can be shortened considerably in this case because large constants are super solutions of the equation. \\

\noindent {\bf Case c: $\Omega_0 = \Omega$} 
 
In this case, the solutions $u_{k_i}$ converge to the sought-after solution $u : \Omega \to \R$ by (\ref{item2}) and (\ref{item11}).  

\begin{remark} Let $(M, g)$ be a complete Riemannian manifold of dimension $n$ with $2 \leq n \leq 7$. Let $H \in \C^\infty(M)$ and let $\Omega \subset M$ be a bounded domain whose boundary can be written as the disjoint union of hypersurfaces $\partial_- \Omega$ and $\partial_+ \Omega$ such that $H_{\partial_- \Omega}  (x)> H(x)$ with respect to the unit normal pointing into $\Omega$ and such that $H_{\partial_+ \Omega}(x) < H (x)$ with respect to the unit normal pointing out of $\Omega$. There exists an open subset $U \subset \Omega$ containing a neighborhood of $\partial_- \Omega$ whose boundary in $\Omega$ is a smooth hypersurface $\Sigma$ whose mean curvature at $x \in \Sigma$  with respect to the unit normal pointing out of $U$ equals $H(x)$. Such a set $U$ can be found by minimizing the functional 
\begin{align} \label{eqn:brane} U \mapsto \H^{n-1}_g (\Omega \cap \partial^*U) - \int_U H (x) d \L^n_g (x).\end{align}
This was proven by M. Fuchs in \cite[Theorems 2.1 and 4.1]{Fuchs:1991}. The existence of a hypersurface with prescribed mean curvature $H$ also follows from the non-variational approach in  \cite{Andersson-Metzger:2009, Eichmair:2009-Plateau}, noting that such surfaces are MOTSs in the initial data set $(M, g, k= - H/(n-1) g)$. The proofs in \cite{Andersson-Metzger:2009, Eichmair:2009-Plateau} proceed by constructing a limit of solutions of regularized Jang equations whose boundary values diverge to plus and minus infinity near $\partial_+ \Omega$ and $\partial_- \Omega$ respectively.  The observation that the capillarity term in the regularized Jang equation contributes ``with a good sign" in the flux integrals (\ref{eqn:flux-}) and (\ref{eqn:flux+}) that we exploited in the proof of Theorem \ref{thm:mainjss} can be used to show that the boundaries of prescribed mean curvature arising in this way also minimize (\ref{eqn:brane}). Cf. \cite[Remark 3.2]{Eichmair:2009-Plateau}.
\end{remark}
 

\appendix

\section {Remark on the Interior gradient estimate for the prescribed mean curvature equation} \label{sec:ige}

Let $(M, g)$ be a complete Riemannian manifold, let $\Omega \subset M$ be a non-empty open subset, and let $u \in C^{2, \alpha}_{loc} (\Omega)$ be a solution of the prescribed mean curvature equation 
\begin{align} \label{eqn:app:pmc}
\div \left( \frac{D u}{ \sqrt {1 + |D u|^2}} \right) = H + t \, u
\end{align}
on $\Omega$
where $t \in (0, 1)$ and where $H$ is a bounded, locally Lipschitz function. The geometric operators appearing in this equation are with respect to the metric $g$, so that the left hand side represents the scalar mean curvature of the graph $G(u) = \{(x, u(x)) : x \in \Omega\} \subset \Omega \times \R$ with respect to the product metric $g + d x^{n+1} \otimes dx^{n+1}$ and the downward pointing unit normal. Then, given $x \in \Omega$, the gradient of $u$ at $x$ is bounded in terms of the geometry of $(M, g)$ near $x$, the size of $H$ in the Lipschitz norm near $x$, and a bound on $u$ near $x$. This interior gradient estimate is derived in \cite{Simon:1976-int} using iteration techniques, and e.g. in \cite{Korevaar:1986}, \cite{Korevaar:1987}, \cite{Korevaar-Simon:1988}, \cite{Wang:1998}, \cite{Spruck:2007} using the maximum principle in a clever way. The estimates for the gradient obtained in these references are explicit. \\

Here, we include a short, indirect, and conceptually simple proof of this well-known interior gradient estimate. This proof works when the dimension $n$ of $M$ lies in the range $2 \leq n \leq 7$. It is based on the regularity theory of (almost) minimizing boundaries. The ingredients of the proof are classical. A very similar argument was employed in the proof of Theorem 4.2 in \cite{Giusti:1978}. For convenience, our references below are to the recent paper \cite{Eichmair:2009-Plateau} whose notation and language are compatible with ours here. There, the reader will find further references to the literature. The line of reasoning here should also be compared to the application of the Harnack inequality in the corollary of Theorem 4.2 of \cite{Simon:1977}. \\ 

Suppose the estimate fails (with a counterexample in dimension $2 \leq n \leq 7)$. Then there exist $(M, g)$, $\Omega \subset M$, $x \in \Omega$, and $H$ as in the statement, a precompact open subset $B \subset \Omega$ containing $x$, and a sequence of solutions $u_k \in C^{2, \alpha}_{loc} (\Omega)$ of (\ref{eqn:app:pmc}) with $t = t_k \in (0, 1)$ such that $u_k(x) = 0$, such that $|u_k (y)| \leq T$ for some $T > 0$ and all $y \in B$, and such that $|D u_k(x)| \to \infty$ as $i \to \infty$. Let $G_k$ denote the graph of $u_k$ above $B$, and let $\nu^{n+1}_k = (1 + |D u_k|^2)^{-1/2}$ denote the vertical component of its upward pointing unit normal. By assumption, the mean curvatures of the graphs $G_k$ are uniformly bounded independently of $k$. Hence the graphs $G_k$ are all contained in a precompact class of almost minimizing (relative) boundaries in $B \times \R$ (cf. \cite[Appendix A and Remark 4.1]{Eichmair:2009-Plateau}). Moreover, there exists a constant $\beta >0$ depending only on the Lipschitz norm of $H$ on $B$ such that $\Delta_{G_k} \nu^{n+1}_k \leq \beta \nu^{n+1}_k$ (cf. (\ref{app:Jacobi}) and (\ref{app:decompose}) below) holds weakly on $G_k$. Moreover, $(x, 0) \in G_k$ for every $i$. Let $G$ be a subsequential varifold limit of $G_k$ as $k \to \infty$. Then $G$ is a properly embedded $C^{2, \alpha}$ two-sided hypersurface in $B \times \R$, and the subsequence of the $G_k$ approaches $G$ in $C^{2, \alpha}$ on compact subsets of $B \times \R$. Let $\nu^{n+1}$ denote the vertical component of the ``upward" unit normal of $G$ (the orientation is inherited from $G_k$). Then $\nu^{n+1}$ is non-negative and $\bar \Delta \nu^{n+1} \leq \beta \nu^{n+1}$ weakly on $G$. By the strong maximum principle, on any component of $G$, $\nu^{n+1}$ is either everywhere positive (so that the component is a graph) or everywhere vanishing (so that the component is a vertical cylinder). See \cite[Section 4.2]{Massari-Miranda:1984} and \cite[Lemma 2.3]{Eichmair:2009-Plateau} for references and a more detailed discussion of this geometric Harnack principle. Since $|D u_k (x)| \to \infty$, it follows that the connected component of $G$ containing $(x, 0)$ is cylindrical. This is impossible, since all the graphs $G_k$ and consequently also $G$ are contained in the slab $M \times [-T, T]$. \\

For convenient reference, and to rectify a mistake in \cite{Eichmair:2009-Plateau}, we include the following lemma:        

\begin{lemma}
Let $(M, g)$ be a complete Riemannian manifold, let $\Omega \subset M$ be non-empty and open, let $u \in \C^2 (\Omega)$, and let $G = \{ (x, u(x)) : x \in \Omega\} \subset M \times \R$ be its graph. Let $v = \sqrt {1 + |D u|^2}$ and $\nu = D u / v$. Then 
\begin{align} \label{app:Jacobi} 
\bar \Delta \frac{1}{v} + (|h|_{\bar g}^2 + \Ric (\nu, \nu) + \nu(H (u))) \frac{1}{v} = 0 
\end{align}
holds weakly on $G$ where $\bar \Delta$ is the (non-positive) Laplace-Beltrami operator of $G$ with respect to the metric $\bar g$ induced on $G$ from $(M \times \R, g + d x^{n+1} \otimes dx^{n+1})$, $|h|_{\bar g}$ is the length of its second fundamental form, $H (u)= \div \nu$ is its mean curvature, and $\Ric$ is the Ricci curvature tensor of $(M, g)$.  

Let $F : T M \times \R \to \R$ be a locally Lipschitz function that is non-decreasing in its last argument and such that $H(u)= F (\nu, u)$. There exist a measurable locally bounded function $f$ and a measurable locally bounded vector field $X$ on $G$ such that 
\begin{align} \label{app:decompose} 
\nu (H (u)) \frac{1}{v} \geq - \frac{f}{v} +  \bar g ( X, \overline \nabla \frac{1}{v})
\end{align}
holds almost everywhere on $G$. The bounds for $f$ and $X$ on a compact subset $K \subset \Omega \times \R$ depend only on the projection of $K$ to the base, as well as the local geometry of $(M, g)$ and the Lipschitz norm of $F$ on that set. The vector field $X$ is the tangential part of a horizontal vector field on $M \times \R$ that is independent of the vertical variable and of $u$. 
If $F : T M \times \R \to \R$ does not depend on the fiber variable in $T M$ we can take $X = 0$.
If we let $w = \log v$, then  
\begin{align} \label{app:w}
\bar \Delta w = |\overline \nabla w|_{\bar g}^2 + |h|^2_{\bar g} + \Ric(\nu, \nu) - f - \bar g ( X, \overline \nabla w)
\end{align}
holds weakly on $G$.
\begin{proof} The Jacobi identity (\ref{app:Jacobi}) is derived carefully in \cite[(2.18)]{Schoen-Yau:1981-pmt2}. For the second part of the lemma, choose local coordinates $(x^1, \ldots, x^n)$ on $\Omega$ so that $F = F(x^1, \ldots, x^n, p^1, \ldots, p^n, z)$. Note that $F_{p^\ell} dx^\ell$ is defined independently of coordinates. Following the notation of \cite{Spruck:2007}, we denote the components $v^{-1} u^i$ of $\nu$ by $\nu^i$ (indices raised with respect to $g$). The induced metric on the graph has components $\bar g_{ij} = g_{ij} + u_i u_j$. The components of its inverse are $\bar g^{ij} = g^{ij} - \nu^i \nu^j$. Then 
\begin{align*}
\nu (H (u)) &= \nu F(\nu, u) = \nu^i F_{p^j} (\nu^j)_{;i} +  \nu^i ( F_{x_i} - F_{p_j}   \Gamma_{ im}^j \nu^m ) + v^{-1}|Du|^2 F_z  \\ &\geq \nu^i F_{p^j} (\nu^j)_{;i} +  \nu^i ( F_{x_i} - F_{p_j}   \Gamma_{ im}^j \nu^m )
\end{align*} where covariant derivatives and the Christoffel symbols are with respect to $g$. The second term on the right is defined independently of our choice of coordinate system; it is our function $f$. For the first term, note that 
$$v^{-1} \nu^i (\nu^j)_{;i} F_{p_j}= v^{-2} \nu^i u_{m;i} \bar g^{mj} F_{p_j} = v^{-3} u^i u_{m;i} \bar g^{mj} F_{p_j} = - ( v^{-1})_m \bar g^{mj} F_{p_j}.$$ We let $X^m = -  \bar g^{m j} F_{p_j}$. Note that $|X|_{\bar g}^2 \leq F_{p_i} F_{p_j} g^{ij}$.
\end{proof}
\end{lemma}

\begin{remark} The derivation of \cite[(6)]{Eichmair:2009-Plateau} is flawed, as was pointed out to the first named author by Zhuobin Liang. The correct ``ensuing differential inequality" is given by (\ref{app:decompose}) along with (\ref{app:Jacobi}). (The gradient term is missing in \cite{Eichmair:2009-Plateau}.) This introduces two changes in \cite{Eichmair:2009-Plateau}: In the proof of Lemma 2.1, a gradient term $\bar g (X, \overline \nabla \eta)$ should be added to the line $- \beta \eta + \Delta_\Sigma \eta$. This introduces a harmless additional term of order $K$ in the estimate that follows in the proof of Lemma 2.1 in \cite{Eichmair:2009-Plateau}. Similarly, a gradient term should be added to the differential inequality in the hypotheses of Lemma 2.3 in \cite{Eichmair:2009-Plateau}. Again, this does not affect the proof. The inequality $(6)$ in \cite{Eichmair:2009-Plateau} was carried over as $(2)$ in \cite{Eichmair:2010}. That inequality should be replaced by (\ref{app:w}), which implies that 
$$\overline \div ( \overline \nabla w - \frac{1}{2} X) - |\overline \nabla w - \frac{1}{2} X |^2_{\bar g} - |h|_{\bar g}^2 = \Ric(\nu, \nu) - f - \frac{1}{4} |X|_{\bar g}^2 - \frac{1}{2} \overline \div  X$$ 
holds weakly on $G$. Multiply this inequality by $\phi^2$ where $\phi \in \C_c^1(\Omega \times \R)$, integrate over $G$, and integrate by parts in a standard way (cf. \cite[Proposition 1]{Schoen-Yau:1981-pmt2}, \cite[Proof of Theorem 2.1]{Galloway-Schoen:2006}), to obtain 
\begin{align} \label{eqnapp:stability} \int_{G} |h|_{\bar g}^2 \phi^2 \leq  (1 + \epsilon ) \int_{G} |\overline \nabla \phi|_{\bar g}^2 + C_\epsilon \int_{G} \phi^2. 
\end{align}
Here, $C_\epsilon$ is a constant that only depends on a given $\epsilon >0$, the projection of the support of $\phi$ to the base, and the size of the Ricci tensor of $(M, g)$ and the derivatives of $F$ on that projection. For appropriate choice of $\epsilon >0$, this inequality is equivalent to estimate $(3)$ in \cite{Eichmair:2010}.  
\end{remark}


\section {Domains of solutions of prescribed mean curvature equations with infinite boundary data} \label{sec:boundary}

Let $(M, g)$ be a complete Riemannian manifold of dimension $n$, $2 \leq n \leq 7$, let $\emptyset \neq \Omega \subsetneq M$ be a non-empty open subset, and let $u \in \C^2(\Omega)$ be a solution of the equation 
\begin{align*} \label{eqn:boundary} \div \left( \frac{D u}{\sqrt {1 + |D u|^2}}\right) = H_0 \end{align*}
where $H_0 \in \R$ is a constant. Let $U \subset M$ be an open set such that $U \cap \partial \Omega \neq \emptyset$ and assume that $\lim_{x \to y, x \in \Omega} u(x) = \infty$ for every $y \in U \cap \partial \Omega $, where $\partial \Omega$ is the topological boundary of $\Omega$ in $M$. \\

The topological boundary $U \cap \partial \Omega$ of $\Omega$ relative to $U$ has the structure of a smooth two-sided properly immersed constant mean curvature hypersurface $\iota: \Sigma \to U$. More precisely, given $y \in U \cap \partial \Omega$, there exists an open neighborhood $\mathcal{O}_y$ of $y$ in $U$ and a diffeomorphism $\phi_y: \mathcal{O}_y \to B^{n-1}_1(0) \times (-1, 1)$ with $\phi_y(y) = 0$ such that one of the following conditions hold: 
\begin{enumerate} [(i)]
\item There exists a smooth function $f: B_1^{n-1} (0) \to \R$ with $f(0) = 0$, $df (0) = 0$, and $|d f (x)| <1$ for all $x \in B_1^{n-1} (0)$ such that $\phi_y(\Omega \cap \mathcal{O}_y) = \{ (x, t) \in B_1^{n-1}(0) \times (-1, 1) : t < f(x)\}$.
\item There exist smooth functions $f^1, f^2: B_1^{n-1} (0) \to \R$ with $f^i(0) = 0$, $df^i (0) = 0$, $|d f^i (x)| < 1$, and $f^1(x) \leq f^2(x)$ for all $x \in B_1^{n-1} (0)$ and $i \in \{1, 2\}$ such that $\phi_y(\Omega \cap \mathcal{O}_y) = \{ (x, t) \in  B_1^{n-1}(0) \times (-1, 1) : t < f^1(x) \text{ or } f^2(x) < t\}$.
\end{enumerate}
Moreover, for every $\psi \in \C_c^1(\Sigma)$ we have that 
\begin{eqnarray} \label{eqn:stabilitycross} 
\int_\Sigma (|h|^2 + \Ric_g (\nu, \nu)) \psi^2  \leq \int_\Sigma |\bar D  \psi|_{\bar g}^2.
\end{eqnarray}
Here, $h$ is the second fundamental form and $\nu$ is the unit normal vector field of the immersion $\iota: \Sigma \to U$, integration is with respect to the volume form of the pull-back metric $\bar g = \iota^* g$, and $\bar D$ is the gradient with respect to $\bar g$. \\

For $n=2$, the conclusions here are (roughly) \cite[p. 329]{Jenkins-Serrin:1966} (for $H_0 = 0$) and \cite[Theorems  6.1 and 6.2]{Serrin:1970} (for $H_0 \neq 0$). For $n \geq 2$, it is shown in Chapter 8 of \cite{Spruck:1972} that $U \cap \partial \Omega$ is a hypersurface of constant mean curvature $H_0$ \emph{assuming} that $\partial \Omega$ is a $\C^2$ hypersurface; see also the analysis of ``extremal domains" in \cite{Giusti:1978}. \\

Below we sketch how these conclusions follow from arguments as applied in \cite{Eichmair:2009-Plateau, Eichmair:2010}. \\   

Note that the graph of $u$, $G = \{ (x, u(x)) : x \in U \cap \Omega\}$, is a properly embedded hypersurface in $U \times \R$. Let $V \Subset U$ be a non-empty open set with smooth boundary that is compactly contained in $U$. Using the divergence theorem for the vector field $(1 + |Du|^2)^{-1/2}  (-Du, 1)$ in the region $\{ (x, x^{n+1}) : x \in \Omega \text{ and } u(x) \leq x^{n+1}\} \cap (V \times [S, T])$ for regular values $S, T$ of $u$ with $S < T$, we obtain that
$$\H^{n}_{g + dx^{n+1} \otimes dx^{n+1}} (G \cap (V \times [S, T])) \leq 2 \L^n_{g} (V) + |H_0| (T -S) \H_g^{n-1} (\partial V).$$ In particular, we obtain locally uniform area bounds for $G$ and its vertical translates. \\

Let $\nu^{n+1} = (1 + |Du|^2)^{-1/2}$ denote the vertical component of the (upward) unit normal vector field $\nu = (1+|Du|^2)^{-1/2} (- D u, 1)$ of $G$. The Jacobi identity (expressing the fact that vertical translation of $G$ leaves the mean curvature constant) implies that 
\begin{eqnarray} \label{Jacobiidentity} \bar \Delta (  \nu^{n+1} )+ (|h|^2 + \Ric_{g + dx^{n+1} \otimes dx^{n+1}}(\nu, \nu)) \nu^{n+1} = 0.\end{eqnarray}
Here, $\bar \Delta$ is the (non-positive) Laplace--Beltrami operator of $G$ with respect to its induced metric $\bar g$. Let $\psi \in \C^1_c(U \times \R)$ be a test function, multiply (\ref{Jacobiidentity}) by $\psi^2 (\nu^{n+1})^{-1}$, integrate over $G$, integrate by parts, and use the elementary estimate $- \bar g ( \bar D (\psi^2 (\nu^{n+1})^{-1}), \bar D \nu^{n+1}) \geq - |\bar D \psi^2|_{\bar g}$ to obtain that
\begin{eqnarray} \label{eqn:stabilityG} \int_{G} (|h|^2 +  \Ric_{g + dx^{n+1} \otimes dx^{n+1}}(\nu, \nu)) \psi^2 \leq \int_G |\bar D \psi|_{\bar g}^2.\end{eqnarray}
The stability--based regularity theory of \cite{Schoen-Simon:1981} applies and provides curvature estimates for $G$. In fact, the curvature of $G \cap (V \times [T, T+1])$ is bounded for every $V \Subset U$ independently of $T\in \R$. (The case $n=7$ requires an additional argument, cf. Remark 4.1 in \cite{Eichmair:2009-Plateau}.) Let $G_T = \{(x, u_T(x)) : x \in \Omega\}$ with $u_T = u - T$. Note that these hypersurfaces are naturally ordered and that their area and curvature is bounded independently of $T$ on compact subsets of $U \times \R$. Hence we can pass them to a geometric (varifold) limit as $T \nearrow \infty$. It is elementary to check that the support of this limit equals $(U \cap \partial \Omega) \times \R$. The arguments in  Appendix A of \cite{Eichmair:2010}, in particular Corollary A.1 and Remark A.3, explain carefully how the  asserted structure of the cross-section of this geometric limit follows from this, in particular why no more than two sheets can come together. That the stability  property (\ref{eqn:stabilityG}) passes to the cross-section in the form (\ref{eqn:stabilitycross}) is checked exactly as in \cite[p. 254]{Schoen-Yau:1981-pmt2}. \\

Similar arguments characterize the boundary of open subsets $\Omega$ of a Riemannian manifold $(M, g)$ that support solutions $u: \Omega \to \R$ with infinite boundary values on $U \cap \partial \Omega$ of the prescribed mean curvature equation 
\begin{align*} \div \left( \frac{D u}{\sqrt {1 + |D u|^2}}\right) = F \left( x, \frac{D u }{\sqrt{1 + |D u|^2}} \right). \end{align*} Here, $F$ is a smooth function on $M \times \mathbb{S}^{n-1} (M)$, where $\mathbb{S}^{n-1}(M)$ is the unit sphere bundle of $(M, g)$. The boundary of $\Omega$ in $U$ then has the structure of a two-sided immersion $\iota: \Sigma \to U$ such hat $H_\Sigma(\iota (\sigma)) = F(\iota(\sigma), \nu(\sigma))$, where $\nu$ is the unit normal field along the immersion pointing out of $\Omega$. Instead of (\ref{eqn:stabilitycross}) one obtains that for every $\epsilon >0$ there exists a constant $C_\epsilon = (\epsilon, |\Ric_M|, |H|_{\C^1})$ such that 
\begin{eqnarray} \label{eqn:almoststability} (1 - \epsilon) \int_\Sigma |h|^2 \psi^2  \leq \int_\Sigma |\bar D  \psi|_{\bar g}^2 + C_\epsilon \int_\Sigma \psi^2
\end {eqnarray}
for every $\psi \in \C^1_c(\Sigma)$. See \cite[($3$)]{Eichmair:2010} for details. 
 

\section {Equicontinuous vector fields from solutions of the prescribed mean curvature equation} \label{sec:equicontinuous}

The lemma below follows from the compactness and regularity theory for graphs with bounded mean curvature, see e.g.  \cite{Eichmair:2009-Plateau} for precise statements and references. 

\begin{lemma} \label{lem:equicontinuity} Let $(M, g)$ be a complete $n$-dimensional Riemannian manifold, $2 \leq n \leq 7$, let $\Omega \subset M$ be a non-empty open subset, and let $C\geq0$. The collection of continuously differentiable vector fields 
\[
\left\{ \frac{D u}{\sqrt{1 + |D u|^2}} : u \in \C^2(\Omega) \text{ and } \left| \div \left( \frac{D u}{\sqrt{1 + |Du|^2}} \right) \right| \leq C \right\}
\]
is equicontinuous on compact subsets of $\Omega$. 
\end{lemma}

 
\section{Barriers for the Jang equation near infinity} \label{sec:barriers}
 
In \cite[p. 248]{Schoen-Yau:1981-pmt2}, certain rotationally symmetric barriers for the Jang equation were constructed on large subsets of the asymptotically flat ends of three dimensional initial data sets. The following proposition, which we quote from \cite{Eichmair:2011-jangreduction}, is a straightforward extension of the construction in \cite{Schoen-Yau:1981-pmt2} to higher dimensions. 

\begin{proposition} \label{prop:barriers} Fix $\beta \in (2, n)$. For $\Lambda \geq 1$ define $$b_\Lambda (r) = \Lambda \int_{\frac{r}{\Lambda}}^\infty \frac{ds}{\sqrt{s^{2(\beta -1)} - 1}}$$ on $[\Lambda, \infty)$. Then $b_\Lambda$ is continuous and positive, it is smooth on $(\Lambda, \infty)$, and we have that $\frac{d b_\Lambda}{dr} (r)$ tends to $-\infty$ as $r \searrow \Lambda$. There is a constant $c = c(\beta) \geq 1$ such that $b_\Lambda (r) \leq c \Lambda (r/\Lambda)^{2- \beta}$ and such that $b_\Lambda (\Lambda) \geq c^{-1} \Lambda$. 

Let $(M, g, k)$ be an asymptotically flat initial data set of dimension $n$, $n \geq 3$, such that $g^{ij} k_{ij} = O(|x|^{-\beta})$ as $x \to \infty$ in the asymptotically flat ends.\footnote{When $n > 3$, our definition of asymptotic flatness in the introduction shows that there exists some $\beta \in (2, n)$ such that $g^{ij} k_{ij} = O(|x|^{-\beta})$. When $n=3$, this is a genuine additional hypothesis, cf. \cite[(1.4)]{Schoen-Yau:1981-pmt2}.}There is $\Lambda_0 = \Lambda_0 (M, g, k, \beta) \geq 1$ such that for every $\Lambda \geq \Lambda_0$ we have that $H(b_\Lambda (|x|)) + \tr (k) (b_\Lambda(|x|)) > 0$ and $H( - b_\Lambda (|x|)) + \tr (k) (- b_\Lambda(|x|)) < 0$ on $\{x \in M : |x| > \Lambda\}$. 
\end{proposition}

In this paper, we repeatedly use the barriers from Proposition \ref{prop:barriers} along with the maximum principle, exactly as in \cite[p. 249]{Schoen-Yau:1981-pmt2}, to conclude that a solution $u$ of the regularized Jang equation $H (u) + \tr(k)(u) = t \, u$ (with $t \geq 0$) that is defined on the complement of a compact subset of $(M, g, k)$ and which tends to zero at infinity actually lies between $- b_\Lambda (|x|)$ and $b_\Lambda(|x|)$ provided that $\Lambda \geq \Lambda_0$ is large enough so that $\{ x \in M : |x| > \Lambda\}$ is contained in the domain of $u$. 


\bibliographystyle{amsplain}
\bibliography{referencespotpourri}
\end{document}